%% file: jumpsAndCoalescenceFFF.tex
\documentclass[reqno,11pt]{amsart}
\usepackage{amsmath,amssymb,mathtools,dsfont,graphicx}
\usepackage[a4paper, margin = 3cm]{geometry}
\newtheorem{theorem}{Theorem}[section]

\theoremstyle{definition}

\theoremstyle{remark} \theoremstyle{remark}

\numberwithin{equation}{section}

\allowdisplaybreaks[3]

\newcommand{\Rd}{\mathds{R}^d}
\newcommand{\intRd}{\int\limits_{\mathds{R}^d}}

\newcommand{\suml}[1]{\sum\limits_{#1}}
\newcommand{\prodl}[1]{\prod\limits_{#1}}

\newcommand{\bs}{\backslash}
\newcommand{\cOne}{\langle c_1 \rangle}
\newcommand{\cTwo}{\langle c_2 \rangle}
\newcommand{\coMax}{c_1^{\text{max}}}

\newcommand{\phinorm}{\langle \phi \rangle}

\newcommand{\Pexp}{\mathcal{P}_{\exp}}

\title[Jumps and Coalescence in the Continuum]{Jumps and Coalescence in the Continuum: a Numerical
Study}

\author{Yuri Kozitsky}
\address{Instytut Matematyki, Uniwersytet Marii Curie-Sk{\l}odowskiej, 20-031 Lublin, Poland}
\email{jkozi@hektor.umcs.lublin.pl}

\author{Igor Omelyan$^\dag$}\thanks{After this work was submitted Igor Omelyan untimely passed away}
\address{Institute for Condensed Matter Physics, National Academy of Sciences of Ukraine, UA-79011 Lviv, Ukraine}
\email{omelyan@icmp.lviv.ua}

\author{Krzysztof Pilorz}
\address{Instytut Matematyki, Uniwersytet Marii Curie-Sk{\l}odowskiej, 20-031 Lublin, Poland}
\email{krzysztof.pilorz@poczta.umcs.lublin.pl}

\begin{document}


\begin{abstract}

The dynamics is studied of an infinite continuum system of jumping
and coalescing point particles. In the course of jumps, the
particles repel each other whereas their coalescence is free. As the
equation of motion we take a kinetic equation, derived  by a scaling
procedure from the microscopic Fokker-Planck equation corresponding
to this kind of motion. The result of the paper is the numerical
study (by the Runge-Kutta method) of the solutions of the kinetic
equation revealing a number of interesting peculiarities of the
dynamics and clarifying the particular role of the jumps and the
coalescence in the system's evolution. Possible nontrivial
stationary states are also found and analyzed.
\end{abstract}
\maketitle

\section{Introduction}

In a broader sense, a typical kinetic equation  is a nonlinear
integro-differential equation describing the temporal evolution of
the density function of a large (infinite) system of `particles'. At
this level of description, the individual particles are not taken
into account and the system is considered as a medium, entirely
characterized by its aggregate parameters like density. A prototype
example is the celebrated Boltzmann equation \cite{Boltz} devised by
Ludwig Boltzmann in 1872 to describe large systems of physical
particles. Since then this approach has received various
applications ranging from the theory of multiple-lane vehicular
traffic \cite{Prig} to the description of evolving ecological
systems \cite{Adams,BP1,Mu,Neuhauser}. Usually, kinetic equations
are devised with the help of phenomenological or heuristic
arguments, and thus are only loosely related to so called `first
principles', e.g., by taking into account appropriate conservation
laws and symmetries. Due to Bogoliubov's pioneering works
\cite{Bogol}, see also \cite{Cer}, it has become clear that the
Boltzmann equation can be derived (by a certain decoupling or
truncating procedure) from an infinite chain of linear equations --
the BBGKY hierarchy -- that describes the microscopic evolution of a
particle system, see e.g., \cite{Uchiyama}. This approach was then
extended to deriving kinetic equations describing ecological systems
from the corresponding microscopic equations of their random
evolution \cite{FKKK,Mu}.

\begin{figure}[h]
\centerline{\includegraphics[width=0.5\textwidth]{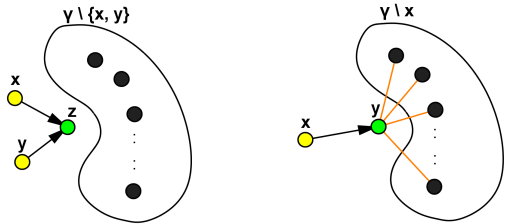}}
\caption{Elementary acts of evolution}\label{rys20}
\end{figure}

In view of various applications -- also but not only those mentioned
above -- there exists a permanent interest to the evolution of
statistically large systems in the course of which the constituents
can merge. As an example, in ecological models merging can be used to 
describe the predation \cite{Delius}.
The Arratia flow \cite{Arrat} provides an example of the
motion of this sort. Its recent study can be found in
\cite{Beres,Kovn,KovnR,LeJ} and in the works quoted therein. In
Arratia's model, an infinite number of Brownian particles move in
$\mathds{R}$ independently up to their collision, then merge and
move together as single particles. Correspondingly, the description
of this motion is performed in terms of diffusion processes. In an
accompanying work \cite{KoP}, we propose an alternative model of
this kind. In this model, an infinite system of point particles
located in $\mathds{R}^d$, $d\geq 1$ undergo random evolution
consisting in the following two elementary acts, see Fig.
\ref{rys20}:
\begin{itemize}
  \item[(a)] Two particles (located at $x$ and $y$) merge into a particle
(located at $z$) with intensity (probability per time) $c_1(x,y;z)$
-- independent of the remaining particles. Thereafter, this new
particle participates in the motion.
\item[(b)] Similarly as in the Kawasaki model \cite{Kawasaki2}, single particles
perform random jumps with repulsion acting on the target point.
\end{itemize}
For the microscopic description of this model, as the phase space
one employs  the set $\Gamma$ of all locally finite configurations
$\gamma \subset \mathds{R}^d$, see \cite{Kawasaki2,
SpatialEcologicalModel} and the next section. The microscopic states
of the system are then probability measures on $\Gamma$ the set of
which is denoted by $\mathcal{P}(\Gamma)$. The evolution of states
$\mu_0\to \mu_t$ is obtained by solving the Fokker-Planck equation.
The main result of \cite{KoP} is the proof of the existence of the
evolution of this type for a bounded time horizon. However, by
virtue of this result the most important and interesting details of
the collective motion of the system remain unrevealed. The aim of
the present work is to study the mentioned model numerically by
employing the corresponding kinetic equation derived from the
microscopic theory developed in \cite{KoP}. The main questions we
address here are: (a) which peculiarities of the motion are related
to each of the mentioned elementary acts of the evolution; (b) what
is the role of the interaction (repulsion) in the possible
appearance of a spatial heterogeneity in the system. In a sense,
this our research is a continuation of the study in \cite{Omel} --
by similar numerical methods -- of the spatial ecological model the
existential problems of which were settled in
\cite{SpatialEcologicalModel}.

\section{Theoretical Background}

As mentioned above, the microscopic theory of our model is based on
the Fokker-Planck equation. The mesoscopic description employs a
kinetic equation obtained from the corresponding microscopic
evolution equations by a scaling procedure, cf. \cite{Ba,FKKK}. Its
solutions are evolving particle densities that will be the objects
of our numerical study.

\subsection{Microscopic description}

The phase space of the dynamics which we study is the set of locally
finite subsets of $\mathds{R}^d$ -- configurations -- defined as
follows
\[
    \Gamma = \{\gamma \subset \Rd: |\Lambda \cap \gamma| < \infty \text{ for any compact } \Lambda \subset \Rd\},
\]
where $|\cdot|$ denotes cardinality. It is equipped with the vague
(weak-hash) topology  see e.g., \cite{Kawasaki2,FKKK}) and the
corresponding Borel $\sigma$-field $\mathcal{B}(\Gamma)$. This
allows one to employ probability measures $\mu$ defined on $(\Gamma,
\mathcal{B}(\Gamma))$ as states of the considered system. The set of
all such measures is $\mathcal{P}(\Gamma)$. The evolution of the
model is described by the Fokker-Planck equation
\begin{equation}\label{KE}
  \mu_t (F) = \mu_0 (F) + \int_0^t \mu_s (LF) d s,
\end{equation}
in which $F:\Gamma\to \mathds{R}$ is an appropriate test function,
$\mu(F):= \int_{\Gamma} F d \mu$, and the operator $L$ specifies the
model. In our case, it is
\begin{gather}\label{OperatorL}
    (LF)(\gamma) = \suml{\{x,y\} \subset \gamma} \  \intRd c_1(x,y;z) \Big(F\big(\gamma \bs \{x,y\} \cup z \big) - F(\gamma) \Big) dz
    \\[.2cm]
    + \suml{x \in \gamma} \  \intRd \tilde{c}_2(x,y;\gamma) \Big(F\big(\gamma \bs x \cup y\big) - F(\gamma) \Big) dy. \nonumber
\end{gather}
Here $c_1\geq 0$ is the intensity of the coalescence of the
particles located at $x$ and $y$ into a new particle located at $z$.
Note that $c_1$ does not depend on the elements of $\gamma$ other
than $x$ and $y$. For simplicity, we assume that $c_1(x,y;z) =
c_1(y,x;z)= c_1 (x+u, y+u; z+u)$ for all $u\in \mathds{R}^d$, i.e.,
$c_1$ is translation invariant. For more general versions of this
model -- that describe also coalescence with interactions --  see
\cite{KP}. The second summand in (\ref{OperatorL}) describes jumps
performed by the particles. Similarly as in \cite{Kawasaki2}, we
take it in the form
\[
    \tilde{c}_2(x,y;\gamma) = c_2(x-y) \prodl{u \in \gamma \bs x} e^{-\phi(y-u)} ,
\]
with $\phi$ and $c_2$ being the repulsion potential and the jump
kernel, respectively. By these assumptions the model is translation
invariant. The functions $c_1$, $c_2$ and $\phi$ take non-negative
values and are supposed to satisfy the following conditions:
\begin{gather} \label{ModelAssumptions}
    \int\limits_{(\mathds{R}^d)^2} c_1(x_1,x_2;x_3) dx_i dx_j = \cOne \ < \infty,
    \\[.2cm]
    \coMax := \sup_{x, y \in \Rd} \intRd c_1(x,y;z) dz < \infty \nonumber,
    \\[.2cm]
    \cTwo \ :=\intRd  c_2(x) dx   < \infty , \quad
    \phinorm := \intRd \phi(x) dx < \infty \nonumber,
    \\[.2cm]
    |\phi| := \sup_{x \in \Rd} \phi(x) < \infty \nonumber.
\end{gather}
As the operator (\ref{OperatorL}) is quite complex, the direct study
of the Fokker-Planck equation (\ref{KE}) is rather inaccessible.
Instead, in \cite{KoP} we realized the following construction. For
$t<T$ (with some $T<\infty$) and $\mu_0$ belonging to a certain
subset of $\mathcal{P}(\Gamma)$, the evolution $\mu_0 \to \mu_t$ was
obtained as the evolution $k_0 \to k_t$ of the correlation functions
corresponding to these states. The basic aspects of this
construction can be outlined as follows. Let $\varOmega$ stand for
the set of all compactly supported continuous functions $\omega
:\mathds{R}^d\to (-1,0]$. Set
\begin{equation*}
  F^\omega (\gamma) = \prod_{x\in \gamma} (1+\omega (x)), \qquad \omega \in \varOmega.
\end{equation*}
Each $F^\omega$ is bounded and continuous, hence integrable for each
$\mu$. Moreover, the collection $\{F^\omega:\omega \in \varOmega\}$
is a measure-defining class, cf. \cite[page 79]{Dawson}. The set of
measures $\mathcal{P}_{\rm exp}\subset \mathcal{P}(\Gamma)$ we will
work with is defined by the condition that its members enjoy the
following property: the map $\varOmega \ni \omega \mapsto
\mu(F^\omega)\in \mathds{R}$ can be continued to an exponential type
entire function defined on $L^1(\mathds{R}^d)$. Then, for $\mu\in
\mathcal{P}_{\rm exp}$, we set $B_\mu(\omega) = \mu (F^\omega)$ and
derive $\widetilde{L}$ from $L$ according to the rule
$(\widetilde{L}B_\mu)(\omega) = \mu(LF^\omega)$. Thereafter, we
construct the evolution $B_{\mu_0}\to B_t$ by solving the
corresponding evolution equation. The next (and the hardest) part of
this scheme is to prove that $B_t = B_{\mu_t}$ for a unique $\mu_t
\in \mathcal{P}_{\rm exp}$. The advantage of using $\mathcal{P}_{\rm
exp}$ is that, for each of its members, the function $B_\mu$ admits
the representation
\begin{eqnarray}
  \label{o4}
B_\mu (\omega) & = & 1 + \sum_{n=1}^\infty \frac{1}{n!}
\int_{(\mathds{R}^d)^n} k^{(n)}_\mu (x_1, \dots , x_n) \omega(x_1)
\cdots \omega(x_n)d x_1 \cdots d x_n .
\end{eqnarray}
Here $k^{(n)}_\mu$ is the $n$-th order correlation function of state
$\mu$. It satisfies the Ruelle bound \cite{Ruelle}
\[
0 \leq k^{(n)}_\mu (x_1, \dots , x_n) \leq \varkappa^n,
\]
with an appropriate $\varkappa>0$. Since each $k^{(n)}$ is defined
by (\ref{o4}) only Lebesgue-almost everywhere, the latter estimate
yields $k^{(n)}\in L^\infty ((\mathds{R}^d)^n)$, and $k^{(n)}$,
$n\geq 2$ is symmetric with respect to the interchange of $x_i$. For
a compact $\Lambda \subset \mathds{R}^d$,
\[
\mu(N_\Lambda)= \int_{\Gamma} N_\Lambda (\gamma) \mu (d \gamma) =
\int_{\Lambda} k_\mu^{(1)} (x ) d x
\]
is the expected value of the number of points contained in $\Lambda$
if the system is in state $\mu$. Here $N_\Lambda(\gamma)
=|\gamma\cap\Lambda|$ is the number of the elements of $\gamma$
contained in $\Lambda$. That is, $k_\mu^{(1)}$ is the particle
density in state $\mu$. Note that $\mu(N_\Lambda)$ may be infinite
for a non-compact $\Lambda$, which would indicate that the system is
infinite in state $\mu$. By the estimate above we have that
\begin{equation}
  \label{o5}
  \|k^{(n)}_\mu\|_{L^\infty
((\mathds{R}^d)^n)} \leq \varkappa^n, \qquad n \in \mathds{N}.
\end{equation}
Let $\Gamma_0$ stand for the set of all finite configurations. It is
a measurable subset of $\Gamma$, equipped with the topology induced
thereon by the vague topology of $\Gamma$. The elements of
$\Gamma_0$ will usually be denoted by $\eta$. Thus, $\eta = \{x_1 ,
\dots , x_n\}$ with distinct $x_j$, and $n=|\eta|$ is the number of
points in $\eta$. Let $k_\mu : \Gamma_0 \to \mathds{R}$ be defined
by $k_\mu (\eta) = k^{(n)}_\mu (x_1, \dots , x_n)$ for $\eta$ as
above. This is the correlation function of state
$\mu\in\mathcal{P}_{\rm exp}$ which characterizes it in a complete
way. For instance, for a Poisson measure $\pi_\rho$ with density
$\rho:\mathds{R}^d\to [0,+\infty)$ we have
\begin{equation}
  \label{CPoi}
k_{\pi_\rho}(\eta) = \prod_{x\in \eta} \rho(x).
\end{equation}
That is, $\pi_\rho$ is completely characterized by its density
(intensity function) $\rho$, cf. \cite[page 45]{Dawson}. By
(\ref{o5}) we conclude that $k_\mu \in \mathcal{K}_\vartheta$ with
$\vartheta = \log \varkappa$ and $\mathcal{K}_\vartheta$ being a
Banach space of such maps equipped with the norm
\[
\|k\|_\vartheta := \sup_{n\geq 0} \bigg{(} e^{- n \vartheta}
\|k^{(n)}_\mu\|_{L^\infty ((\mathds{R}^d)^n)}\bigg{)}.
\]
Then the states $\mu_t \in \mathcal{P}_{\rm exp}$ satisfy the 
Fokker-Planck equation \eqref{KE} with
$L$ given in (\ref{OperatorL}) and $F=F^\omega$ if their correlation
functions $k_t$ satisfy
\begin{equation}
  \label{CFE}
  \frac{d}{dt} k_t = L^\Delta k_t, \qquad k_t|_{t=0}= k_{\mu_0},
\end{equation}
which, in fact, is an infinite chain of equations for $k_t^{(n)}, n
\in \mathds{N}$. In \eqref{CFE}, $L^\Delta$ has the form, cf.
\cite{KP},
\begin{equation*}
    L^\Delta = L_1^\Delta + L_2^\Delta.
\end{equation*}
Here $L_1^\Delta = L_{11}^\Delta + L_{12}^\Delta + L_{13}^\Delta +
L_{14}^\Delta$ is the part responsible for the coalescence whereas
$L_2^\Delta = L_{21}^\Delta + L_{22}^\Delta$ describes the jumps.
Their summands are:
\begin{gather*}
    (L^{\Delta}_{11} k) (\eta) = \ \ \frac{1}{2} \int\limits_{(\mathds{R}^d)^2} \sum\limits_{z \in \eta} c_1(x,y;z) k(\eta \backslash z \cup \{x,y\})  dx dy,
    \\[.2cm] \nonumber
    (L^{\Delta}_{12} k) (\eta) = - \frac{1}{2} \int\limits_{(\mathds{R}^d)^2} \sum\limits_{x \in \eta} c_1(x,y;z) k(\eta \cup y)  dy dz,
    \\[.2cm] \nonumber
    (L^{\Delta}_{13} k )(\eta) = - \frac{1}{2} \int\limits_{(\mathds{R}^d)^2} \sum\limits_{y \in \eta} c_1(x,y;z) k(\eta \cup x)  dx dz,
    \\[.2cm] \nonumber
    (L^{\Delta}_{14} k )(\eta) =  - \Psi (\eta) k(\eta), \qquad \Psi(\eta):= \int\limits_{\mathds{R}^d}
    \sum\limits_{\{x,y\} \subset \eta} c_1(x,y;z) dz,
\end{gather*}
and
\begin{align*}
    L^{\Delta}_{21} k(\eta) &= \intRd \sum\limits_{y \in \eta}  c_2(x-y) \prodl{u \in \eta \backslash y} e^{-\phi(y-u)} \ (Q_y k)(\eta \backslash y \cup x) dx,
    \\[.2cm]
    L^{\Delta}_{22} k(\eta) &= - \intRd \sum\limits_{x \in \eta} c_2(x-y) \prodl{u \in \eta \backslash x} e^{-\phi(y-u)} \ (Q_y k)(\eta) dy,
\end{align*}
see \cite{KoP} for more detail. In order to figure out the real
meaning of (\ref{CFE}), let us write down the first two members of
this chain of equations. The first one reads
\begin{eqnarray}
  \label{First}
  \frac{d }{dt} k_t^{(1)} (z) & = & \frac{1}{2}
  \int_{(\mathds{R}^d)^2}c_1 (x,y;z) k_t^{(2)}(x,y) d x d y \\[.2cm]
  \nonumber & - &  \int_{(\mathds{R}^d)^2}c_1 (z,y;x) k_t^{(2)}(z,y) d x d y
  \\[.2cm] \nonumber & + & \int_{\mathds{R}^d} c_2 (x-z) \left[(Q_z
  k_t)^{(1)} (x) - (Q_x k_t^{(1)}) (z)\right] d x,
\end{eqnarray}
where
\begin{eqnarray}
  \label{FirstQ}
  (Q_x k_t)^{(1)}(y)& = & k_t^{(1)} (y) + \sum_{n=1}^\infty \frac{1}{n!}
  \int_{(\mathds{R}^d)^n} k_t^{(n+1)} (x_1, \dots , x_n , y) \\[.2cm] \nonumber &\times &  \left(
 \prod_{j=1}^n \left[ e^{-\phi(x- x_j)} -1 \right]\right) d x_1
 \cdots d x_n.
\end{eqnarray}
The second member of the chain is
\begin{eqnarray}
  \label{Second}
\frac{d }{dt} k_t^{(2)} (x,y) & = & \frac{1}{2}  \int_{(\mathds{R}^d)^2}\left[c_1 (u,v;x) k_t^{(3)}(y,u,v) + c_1 (u,v;y)
 k^{(3)}_t (x,u,v)\right] d u d v \qquad \qquad \\[.2cm] \nonumber & - &
 \int_{(\mathds{R}^d)^2}\left[ c_1 (x,u;v) k_t^{(3)} (x,y,u)  + c_1 (y,u;v)  k_t^{(3)} (x,y,u) \right] d u d v \\[.2cm] \nonumber &  - & \Psi (x,y) k^{(2)}_t (x,y)\\[.2cm] \nonumber & + &
\int_{\mathds{R}^d} e^{-\phi(x-y)} \left( c_2 (u-x) (Q_x
k_t)^{(2)}(y,u) + c_2 (u-y) (Q_y k_t)^{(2)}(x,u) \right) du \\[.2cm] \nonumber &
- & \int_{\mathds{R}^d}\left( c_2 (x-u) e^{-\phi(y-u)} +  c_2 (y-u)
e^{-\phi(x-u)} \right)(Q_u k_t)^{(2)}(x,y) d u,
\end{eqnarray}
where
\begin{equation*}
  \Psi(x,y) = \int_{\mathds{R}^d}\left( c_1 (x,y;z) + c_1(y,x;z)
  \right) dz,
\end{equation*}
and
\begin{eqnarray*}
 (Q_z k_t)^{(2)}(u,v) & = & k^{(2)}_t (u,v) + \sum_{n=1}^{\infty}
 \frac{1}{n!}  \int_{(\mathds{R}^d)^n} k_t^{(n+2)} (u,v, x_1 , \dots
 , x_n) \\[.2cm] \nonumber & \times & \left(\prod_{i=1}^n \left[ e^{-\phi(z-x_i)}-1\right]
 \right) d x_1 \cdots d x_n.
\end{eqnarray*}
Noteworthy, unlike to the most of such equations, cf. \cite[eqs. (1)
and (2)]{Omel}, the right-hand sides of (\ref{First}) and
(\ref{Second}) contain correlation functions of all orders.

Let us now return to studying (\ref{CFE}). By the very definition of
the norms $\|\cdot \|_\vartheta$ the Banach spaces
$\mathcal{K}_\vartheta, \vartheta \in \mathds{R}$ constitute an
ascending scale of spaces such that $\mathcal{K}_{\vartheta'}
\hookrightarrow \mathcal{K}_\vartheta$ whenever $\vartheta'<
\vartheta$. Here $\hookrightarrow$ denotes continuous embedding. By
mean of the estimates
\begin{eqnarray*}
 \left|(L^\Delta_{1i}k)(\eta)\right|& \leq & \left(\frac{1}{2} \langle c_1 \rangle e^\vartheta \|k\|_\vartheta \right)
  |\eta| e^{\vartheta |\eta|}, \quad i=1,2,3, \\[.2cm]
 \left|(L^\Delta_{14}k)(\eta)\right| &\leq & \left(\frac{1}{2}  c_1^{\rm max}  \|k\|_\vartheta \right) |\eta|(|\eta|-1)e^{\vartheta
 |\eta|},
\end{eqnarray*}
as well as, cf. \cite[eq. (3.18)]{Kawasaki2},
\begin{equation*}
\left|(L^\Delta_{2}k)(\eta)\right| \leq \left(2 \langle c_2 \rangle
\exp\left(\langle \phi \rangle e^\vartheta \right)\|k\|_\vartheta
\right)|\eta| e^{\vartheta |\eta|},
\end{equation*}
one defines bounded linear operators $L^\Delta:
\mathcal{K}_{\vartheta'} \to \mathcal{K}_\vartheta$ and then places
the Cauchy problem (\ref{CFE}) into the mentioned scale of Banach
spaces. The results of \cite{KoP} are contained in the following
two statements.
\begin{theorem}\label{Th1}
For each $\vartheta_0 \in \mathds{R}$ and $\vartheta_* >
\vartheta_0$, and for an arbitrary $k_0 \in
\mathcal{K}_{\vartheta_0}$, the problem in \eqref{CFE} has a unique
classical solution $k_t \in \mathcal{K}_{\vartheta_*}$ on $[0, T)$
with the bound $T=T(\vartheta_*, \vartheta_0)$ dependent on
 $\vartheta_0$ and $\vartheta_*$.
\end{theorem}
A priori the solution $k_t$ described in Theorem \ref{Th1} need not
be a correlation function of any state, which means that the result
stated therein has no direct relation to the evolution of states of
the system considered. The next statement removes this drawback.
\begin{theorem}\label{Th2}
Let $\mu_0 \in \Pexp$ be such that $k_{\mu_0}\in
\mathcal{K}_{\vartheta_0}$ and $T(\vartheta_*, \vartheta_0)$ be as
in Theorem \ref{Th1}. Then the evolution $k_{\mu_0}\to k_t$
described in Theorem \ref{Th1} has the property: for each $t<
T(\vartheta_*, \vartheta_0)/2$, $k_t$ is the correlation function of
a unique state $\mu_t \in \Pexp$.
\end{theorem}
The fact that the evolution  described in Theorem \ref{Th2} is only
local in time stems from technical limitations of the method used in
the proof of this statement. We believe that this drawback could be
overcome.

\subsection{Mesoscopic description}

Theorem \ref{Th2} states the existence of the evolution $\mu_0 \to
\mu_t$ of the micro-states of our system obtained by solving the
Fokker-Plank equation (\ref{KE}). In order to get more detailed
information regarding this evolution, we pass to the mesoscopic
level based on the kinetic equation which we derive now. Its naive
(one may say \emph{more direct}) version can be outlined as follows.
As is known, Poisson measures (with correlation functions given in
(\ref{CPoi})) correspond to the states of systems of particles
independently distributed over $\mathds{R}^d$. Possible dependencies
in a state $\mu$ can be captured by the deviations of
$k_\mu^{(n)}$'s from the corresponding products of $k_\mu^{(1)}$'s.
In particular, by the truncated second-order correlation function,
cf. \cite{Reb},
\[
\tau^{(2)}_\mu (x,y) = k^{(2)}_\mu (x,y) -  k^{(1)}_\mu (x)
k^{(1)}_\mu (y),
\]
Thus, a naive truncation (called also moment closure \cite{Mu})
consists in putting by force
\[
k_t^{(n)} (x_1 , \dots , x_n) =  k^{(1)}_t(x_1)k^{(1)}_t(x_2) \cdots
k^{(1)}_t(x_n), \qquad n\geq 2,
\]
in the right-hand sides of (\ref{First}), (\ref{FirstQ}), and
forgetting of (\ref{Second}) and of the equations for $k_t^{(n)}$,
$n\geq 3$. In physical language, it corresponds to a so called
``mean-field approximation", cf. \cite{Omel} where we get out beyond
it. After making this \emph{ansatz}, one obtains a single nonlinear
equation, see (\ref{KEa}) below with $\rho_t = k^{(1)}_t$ -- instead
of the infinite chain of linear equations encrypted in (\ref{CFE}).
A more sophisticated version of passing to the kinetic equation in
(\ref{KEa}) is based on the so-called Poisson approximation of the
states, in which each $\mu_t$ is approximated in a certain way by
the Poisson state $\pi_{\rho_t}$ with density $\rho_t$. Namely, we
say that $\mu\in\mathcal{P}_{\rm exp}$ is approximable by $\pi_\rho$
(Poisson-approximable), if for some $\vartheta \in \mathds{R}$ there
exists a continuous mapping $A: [0,1] \rightarrow
\mathcal{K}_\vartheta$ such that $A(0) = k_{\pi_\rho}$ and $A(1) =
k_\mu$. The approximation scheme can be depicted as follows
\begin{center}
    \begin{tabular}{lclcl}
   $\mu_0$ & $\xleftrightarrow[\text{approximation}]{\text{Poisson}}$ &
$\pi_{\rho_0}$ & $\rightarrow$ & $\rho_0$ \\
   $\downarrow$ {\scriptsize FPE} &  &  &  & $\downarrow$ {\scriptsize KE}
\\
   $\mu_t$ & $\xleftrightarrow[\text{approximation}]{\text{Poisson}}$ &
$\pi_{\rho_t}$ & $\leftarrow$ & $\rho_t$ \\
    \end{tabular}
\end{center}
Its precise formulation is given in the next statement.
\begin{theorem}\label{Th3}
    Let $\mu_0 \in \Pexp$ be approximable by $\pi_{\rho_0}$. Then
    $\mu_t$ stated in Theorem \ref{Th2} is approximable by
    $\pi_{\rho_t}$ for $t < T(\vartheta_*, \vartheta_0)/2$, with
    $T(\vartheta_*, \vartheta_0)$ as in Theorem \ref{Th1} and $\rho_t$
    being the solution to kinetic equation (\ref{KEa}) with initial condition $\rho_0$.
\end{theorem}
The proof of this theorem -- quite technical -- is essentially based
on the method developed in \cite{KoP}. It will be presented in a
separate work. Here we only note that the temporal locality in this
statement is also a matter of technical limitations of the method
used in the proof of Theorem \ref{Th2}.

 The kinetic equation corresponding to (\ref{CFE}) is
\begin{eqnarray} \label{KEa}
    \frac{d}{dt} \rho_t (x)& = & \frac{1}{2} \int_{(\mathds{R}^d)^2} c_1(y,z;x) \rho_t (y) \rho_t(z) d y dz  \\[.2cm] \nonumber
    & - & \left(\int_{(\mathds{R}^d)^2} c_1(x,y;z) \rho_t (y)  d y dz \right) \rho_t(x)\\[.2cm] \nonumber
    & + & \int_{\mathds{R}^d} c_2 (x-y) \exp\left(- \int_{\mathds{R}^d} \phi(y-u) \rho_t (u) d u  \right) \rho_y(y) dy \\[.2cm] \nonumber
    & - & \left( \int_{\mathds{R}^d} c_2 (x-y) \exp\left(- \int_{\mathds{R}^d} \phi(y-u) \rho_t (u) d u  \right) dy \right) \rho_t(x).
\end{eqnarray}
Here the first two terms describe the coalescence whereas the
remaining ones correspond to the jumps. This equation is the object
of the numerical study the results of which are presented  in the
next section. Note that a heuristic derivation/justification of this
equation could hardly be done, which once again manifests advances
of our method of deriving kinetic equations from the corresponding
microscopic theories.

\section{Simulations}

\subsection{Notions and techniques}
The numerical study of (\ref{KEa}) is performed by simulations based
on the fourth-order Runge-Kutta method. For further simplicity, we
restrict the study to the one-dimensional case. As to the intensity
functions, we consider the following specific cases. The coalescence
intensity is taken in two forms
\begin{align}
    c_1(x_1, x_2; x_3) &= a_1 (x_1 - x_2) \delta\Big(\frac{x_1 + x_2}{2} - x_3\Big), \label{coalescenceOne}\\
    c_1(x_1, x_2; x_3) &= a_1 (x_1 - x_2) \delta\Big(\ln(e^{x_1} + e^{x_2}\big) - x_3\Big). \label{coalescenceTwo}
\end{align}
In  both cases, $\delta$ stands for the Dirac $\delta$-function,
i.e., $c_1$ is a distribution. In the case of
\eqref{coalescenceOne}, the resulting point of the coalescence is at
the middle of the coalescing particles, which may be a natural
choice for describing by $x,y$ their spatial location. The  case of
\eqref{coalescenceTwo} represents the coalescence with the
conservation of mass in the case where $x= \ln m$, $m$ being
particle's mass. In both cases, instead of the $\delta$-function,
one can consider more smooth functions to take into account possible
dispersion. We adopt these forms to spare the calculation time. Note
that while the case of \eqref{coalescenceOne} can easily be
generalized to $d>1$, the case of \eqref{coalescenceTwo} is
essentially one-dimensional. The kernels $a_1$, $c_2$ and $\phi$
were chosen to be non-negative and symmetric functions -- either
Gaussian \eqref{kernelOne} or simple step-like  \eqref{kernelTwo}:
\begin{align}
    G_{\lambda, \sigma}(x) &= \frac{\lambda}{\sigma \sqrt{2 \pi} } \exp \Big( - \frac{x^2}{2 \sigma^2} \Big) \label{kernelOne},\\
    B_{\lambda, \sigma}(x) &= \frac{\lambda}{2\sigma} I_{[-\sigma, \sigma]}(x)  \label{kernelTwo},
\end{align}
where $I$ stands for indicator, positive $\lambda$ and $\sigma$ are
a strength  and a range parameters, respectively. For both forms of
the kernels, we also use their shifted versions. To introduce them,
we define a symmetry-preserving shift operation $S_h$ \[
    \Big(S_h (f) \Big) (x) = \frac{1}{2} \Big( f(x - h) + f(x + h) \Big), \qquad
    h>0,
\]
and then set
\begin{align*}
    G_{\lambda, \sigma, h} &= S_h \Big(G_{\lambda, \sigma}\Big), \\
    B_{\lambda, \sigma, h} &= S_h \Big(B_{\lambda, \sigma}\Big).
\end{align*}
Note that these choices correspond to the translation invariance of
the system, cf. (\ref{ModelAssumptions}).

In order to imitate the behavior of infinite system, three different
choices of boundary conditions were applied, depending on the
initial state of the system. In the case when in given direction the
entities are absent, the zero Dirichlet boundary condition was used.
In the case of homogeneous distribution in given direction, the
boundary condition was set to be the time-dependent value
corresponding to the analytic homogeneous solution. Finally, for the
cyclic initial condition, the periodic (toroidal) boundary
conditions were applied. In the first two cases an automatic size
adjustment was used to ensure safe distance between the wavefronts
and boundaries. The initial domain of simulations was defined as
segment $[-\frac{L}{2}, \frac{L}{2}]$ with various values of
parameter $L$. In the case of triggering the enlargement mechanism,
the length of the segment was being doubled.

In the following parts, the results of the performed simulations with different
choices of initial conditions and described parameter functions are presented.

\subsection{Jumps in the absence of coalescence}
\begin{figure}[ht]
    \caption{Free jumps without coalescence. Density $\rho_T$ on $[-10, 10]$
        for periodic $B_{1,1}$-type initial condition with period 40. On the left:
        evolution in time with $G_{1,1}$ jump kernel. On the right: comparison at
        $T = 20$ with different choices of jump kernel.}
    \label{FigureFreeJumps1}
    \resizebox{0.45\linewidth}{!}{\input{fig1}}
    \resizebox{0.45\linewidth}{!}{\input{fig2}}
\end{figure}
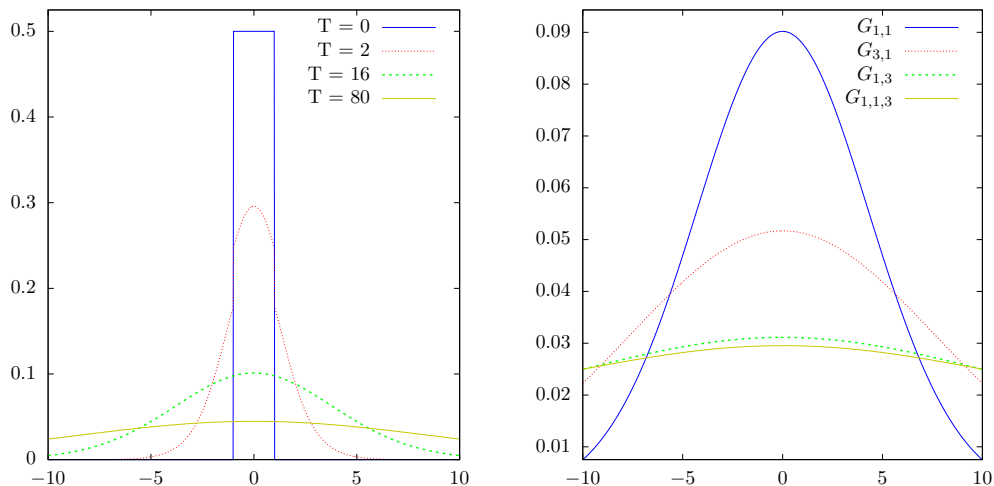

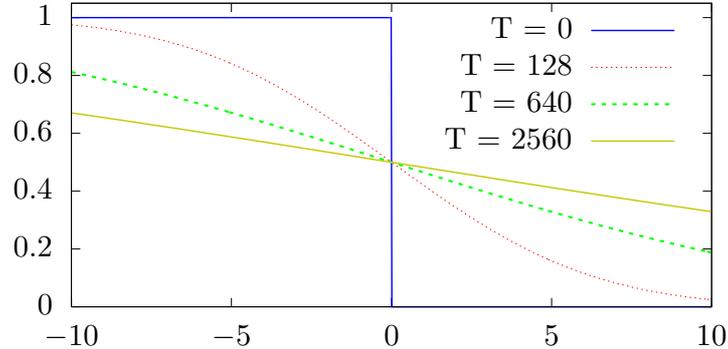
\begin{figure}
    \caption{Free jumps without coalescence. Density $\rho_T$ on $[-10, 10]$ for
        initial condition $I_{(-\infty, 0]}$ with $G_{1,1}$ jump kernel.}
    \label{FigureFreeJumps2}
    \resizebox{0.7\linewidth}{!}{\input{fig3}}
\end{figure}

First, let us analyze the behavior of the system when the entities
are allowed only to jump. In the case of free jumps ($\phi = 0$),
with the exception of some degenerated jump kernels, the system
tends to homogeneity, see Figure \ref{FigureFreeJumps1}, left plot
or Figure \ref{FigureFreeJumps2}.

In the case presented in the Figure \ref{FigureFreeJumps1}, the
initial density function was chosen to be a cyclic step function
$\rho_0(x + 40k) = B_{1,1}(x)$. The discontinuity of the initial
condition dissipates with time, transforming the initial density at
the considered segment $[-20, 20]$ into one resembling a Gaussian,
which tend to more and more homogeneous shape. Increase of the
strength or range parameter, as well as the shift of the jump kernel
result in acceleration of this process, see Figure
\ref{FigureFreeJumps1}, right plot. The toroidal boundary condition
was applied on both sides in order to imitate the cyclic density.

When choosing a different initial condition, for example $I_{(-\infty, 0]}$,
the system is tending to a homogeneity as well, see Figure \ref{FigureFreeJumps2}.
The initial discontinuity vanishes with time and one can observe the rise of
S-shaped densities that preserve at the initial discontinuity point the arithmetic
mean of two initial values (0 to the right and 1 to the left). With the flow of time,
the density function flattens and it is expected that asymptotically the system tends
to this mid-value everywhere. Similarly as in the previous choice of the initial
condition, the higher values of kernel parameters increase the speed of change.
In this case on both sides the Dirichlet boundary condition with auto-enlargement
procedure was used (1 on the left and 0 on the right), as the homogeneous solution
in the absence of coalescence is constant.

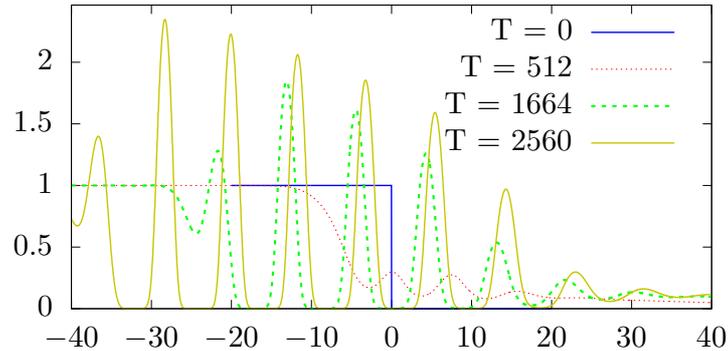
\begin{figure}
    \caption{Jumps with repulsion without coalescence. Density $\rho_T$ on
        $[-40, 40]$ for initial condition $I_{(-\infty, 0]}$ with $G_{1,1,2}$
        jump kernel and $G_{10,1,4}$ repulsion potential.}
    \label{FigureRepulsiveJumps}
    \resizebox{0.7\linewidth}{!}{\input{fig4}}
\end{figure}

The introduction of repulsive effect by setting non-zero $\phi$
usually slows down the flattening process with addition of small
local disturbances in the initial phase. However, in some cases it
may lead to appearance of self-propagating spatial heterogeneity
that seems to drastically change even the asymptotical behavior of
the system. In order to observe such phenomenon, we chose both the
jump kernel $c_2$ and the repulsion potential $\phi$ to be shifted
Gaussians, with the shift $h$ and strength $\lambda$ of repulsion
exceeding the corresponding values for jump kernel. The results are
presented in Figure \ref{FigureRepulsiveJumps}, where the appearance
and propagation of heterogeneity can be observed. One can expect
that the system tends to a non-homogeneous stationary solution. Note
how the initial domain $[-20, 20]$ was automatically enlarged
between $T = 0$ and $T = 512$ in order to avoid artificial behavior
at the boundary and increase simulation accuracy.

\subsection{Pure coalescence}
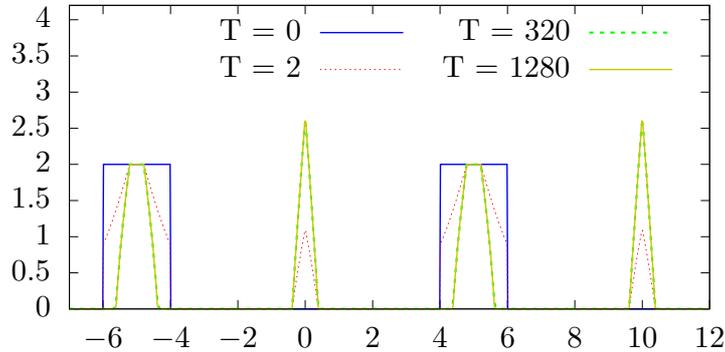
\begin{figure}
    \caption{Pure spatial coalescence. Density $\rho_T$ on [-7, 12] for periodic
        $B_{4,1}$-type initial condition with period 10 and $B_{1,0.8,8}$
        coalescence kernel.}
    \label{FigureCoalescence1}
    \resizebox{0.7\linewidth}{!}{\input{fig6}}
\end{figure}

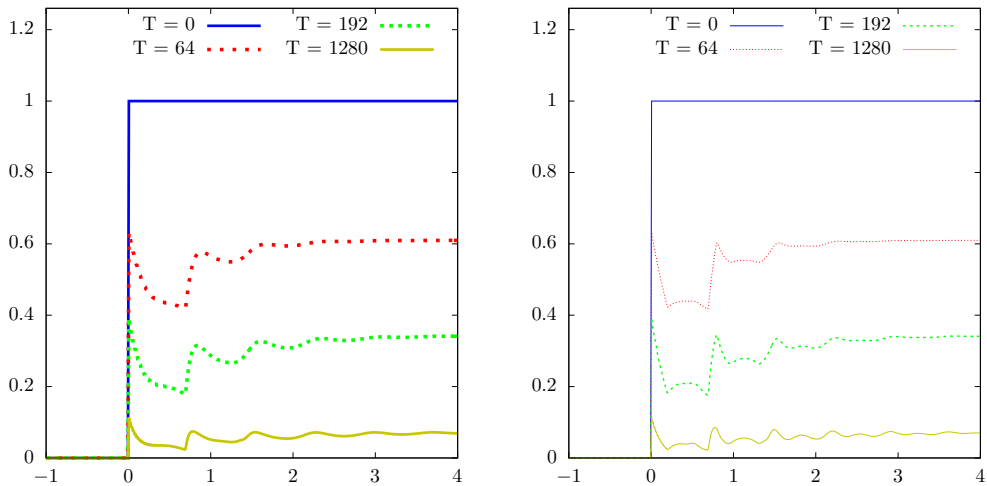
\begin{figure}
    \caption{Pure coalescence with preservation of mass. Density $\rho_T$ on
        $[-1, 4]$ for initial condition $I_{[0,\infty)}$. Coalescence kernel
        used on the left: $G_{0.02,0.2}$, on the right: $B_{0.02,0.2}$.}
    \label{FigureCoalescence2}
    \resizebox{0.45\linewidth}{!}{\input{fig13}}
    \resizebox{0.45\linewidth}{!}{\input{fig7}}
\end{figure}

Let us turn to the pure coalescence without additional jump term. In
the standard cases when coalescence kernel is positive at zero,
independently of the initial condition the system asymptotically
tends to the trivial null density. However, with special choice of
initial condition together with shifted kernel, it is possible to
obtain non-trivial dynamics that seems to be tending to a stationary
state, see Figure \ref{FigureCoalescence1}. Moreover, even when
asymptotic behavior is trivial, the finite time behavior may prove
to be very interesting, see \ref{FigureCoalescence2}.

For the simulation presented in Figure \ref{FigureCoalescence1}, we chose coalescence
intensity of the form \eqref{coalescenceOne} with the shifted step kernel
and the initial density being a periodic step function. The parameters of
coalescence kernel were chosen in such a way that one part of the initial
density is left invariant, while another part diminishes producing peaks between
initial steps. The system seems to approach a non-trivial stationary state,
as the differences between subsequent iterations become very small (compare
density at moments of time $T = 320$ and $T = 1280$). In this case the toroidal
conditions were used at the boundaries of initial domain [-20, 20]. Note that
by changing the shift of coalescence kernel one can easily make the initial
state invariant, e.g. choosing the shift $h = 6$ instead of $h = 8$.

Next, consider the second choice of the coalescence intensity \eqref{coalescenceTwo}.
In this case, the studied property represents the logarithm of mass and during
the action of coalescence no mass is being lost. The simulated density shows
the distribution of mass and how it evolves in time.

As an example, consider initial density $\rho_0$ to be $I_{(-\infty,0]}$ and
observe its change in time, see Figure \ref{FigureCoalescence2}. Pick a weak
coalescence kernel with strength $\lambda = 0.02$ and range $\sigma = 0.2$.
The sudden change from $0$ to $1$ in the density at point $0$ (there are no
entities of mass smaller than 1) produces an interesting irregularity close to
this discontinuity with the flow of time. While the density is approaching zero,
a specific pattern in density remains and propagates with weaker and weaker
amplitude to the right. Note that choice of $B_{0.02,0.2}$ (Figure
\ref{FigureCoalescence2} on the right) kernel produces more irregular
shape than $G_{0.02,0.2}$ (on the left).

Notice that due to the choice of initial condition, despite the mass
preserving choice of coalescence intensity, the density diminishes
to zero pointwise, giving asymptotically at $T \rightarrow \infty$
null density. This behavior shows dishonesty of considered system
with such choice of initial condition. Of course, at any time $T <
\infty$ the total mass is constantly infinite.

\subsection{Jumps and coalescence}
So far we discussed some phenomena that occur separately for jumps
and coalescence. It is interesting that in each considered case, the
introduction of both jump and coalescence nonzero intensities seems
to have a regulating effect on the system that puts it on the path
leading to homogeneity.

\begin{figure}
    \caption{Jumps with shifted repulsion kernel in presence of coalescence.
        Density $\rho_T$ on $-70,10$ for $I_{(-\infty, 0]}$ initial condition.
        Jump kernel $G_{1,1,2}$ with repulsion potential $G_{10,1,4}$ and
        $G_{0.05,1,2}$ coalescence kernel.}
    \label{FigureJumpsAndCoalescence1}
    \resizebox{0.7\linewidth}{!}{\input{fig14}}
\end{figure}
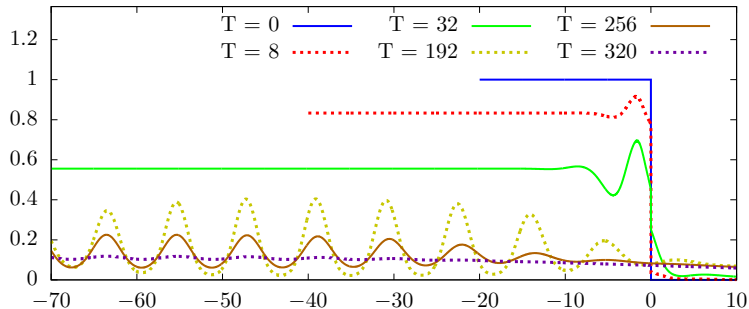

First, consider the case of initial density $\rho_0 = I_{(-\infty,0]}$, where
the choice of shifted repulsion potential significantly stronger than jump kernel
leads to the appearance of periodic spatial heterogeneity (Figure \ref{FigureRepulsiveJumps}).
Introduce to the system relatively weak coalescence kernel $G_{0.05,1,2}$ (with
the choice \eqref{coalescenceOne} of coalescence intensity) in comparison to
$G_{1,1,2}$ jump kernel or $G_{10,1,4}$ repulsion potential. The results are
presented in Figure \ref{FigureJumpsAndCoalescence1}, where the initial simulation
domain was chosen as previously to be $[-20,20]$, Dirichlet zero boundary  condition
was used on the right and time-dependent boundary condition being the solution
to homogeneous problem on the left. Due to the action of coalescence, the density
level reduces and the heterogeneous irregularities start to appear. The addition
of coalescence seems to speed up the process of their formation, cf. density at
$T=192$ in Figure \ref{FigureJumpsAndCoalescence1} and at $T=512$ in Figure
\ref{FigureRepulsiveJumps}. However, under the unceasing influence of coalescence,
the density starts to flatten out (cf. $T=192, 256$ and $320$) and the system obtains
much more homogeneous structure with inevitable null density at the time limit
$T \rightarrow \infty$.

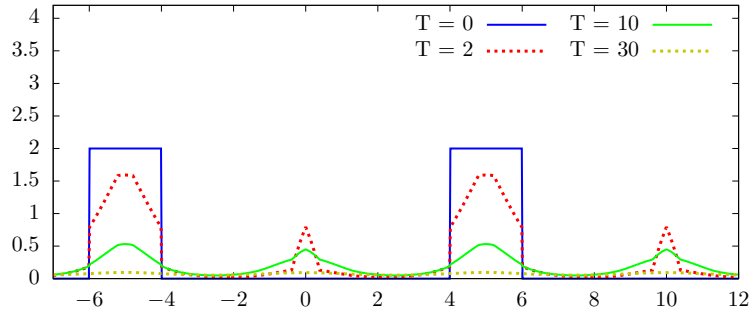
\begin{figure}
    \caption{Coalescence with free jumps. Density $\rho_T$ on $[-7,12]$ for
        $B_{4,1}$-type initial condition with period 10. $B_{1,0.8,8}$ coalescence
        kernel and $J_{0.2,1}$ jump kernel.}
    \label{FigureJumpsAndCoalescence2}
    \resizebox{0.7\linewidth}{!}{\input{fig10}}
\end{figure}

Turn to the case of stationary state appearance for shifted coalescence kernel
in the absence of jumps. The reason of this phenomenon is that after initial dynamic
evolution of density, the shape of coalescence kernel prevents the remaining positive
density areas to reach each other, which results in stagnation of the system, see
Figure \ref{FigureCoalescence1}. The introduction of even small nonzero jump intensity
breaks this impasse. It makes the density in those areas to spread with the flow of time,
which allows the coalescence to operate again, leading the system to null density
limit at $T \rightarrow \infty$, see Figure \ref{FigureJumpsAndCoalescence2}. Note
that the peaks between initial steps still appear, increasing the frequency of higher
density areas for $0 < T < \infty$.

\begin{figure}
    \caption{Coalescence with preservation of mass in presence of free jumps.
        Density $\rho_T$ on $[-1, 4]$ for initial condition $I_{[0,\infty)}$.
        Coalescence kernel $G_{0.02,0.2}$ and $G_{0.01,0.2}$ jump kernel.}
    \label{FigureJumpsAndCoalescence3}
    \resizebox{0.7\linewidth}{!}{\input{fig12}}
\end{figure}
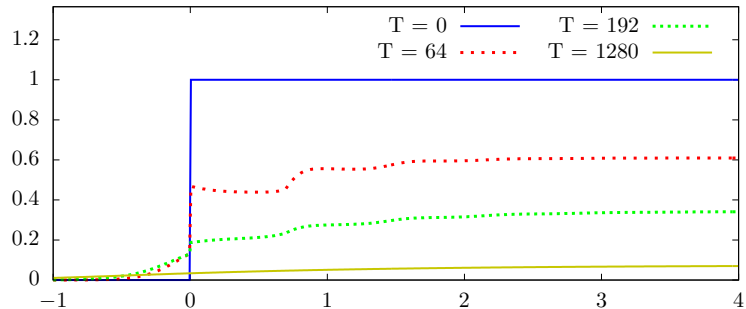

The addition of nonzero jump kernel proves to have a regulating effect in the case
of mass preserving dynamics presented in Figure \ref{FigureCoalescence2} as well.
The irregularities that were appearing and propagating due to the coalescence,
quickly fades in the presence of jumps, see Figure \ref{FigureJumpsAndCoalescence3}.
Introduction of jumps results in more homogeneous tendency to ultimate null density
at time limit. It also leads to appearance of entities with mass smaller than 1.

\section{Conclusions}
The performed simulations show the possibilities of nontrivial dynamics for
the Poisson approximation of coalescing random jumps. The simulated solutions
for the obtained kinetic equation were presented for several cases of initial
conditions and intensities involved. The presented results can be summarized
in the form of following observations:
\begin{enumerate}
    \item In the case of free jumps the system tends to a homogeneous density
    along the whole system. The increase of strength, range or shift of jump
    kernel speeds up this tendency.
    \item For repulsive jumps with shifted repulsion intensity and inhomogeneous
    initial density, persisting spatial heterogeneity may appear that expand with
    the flow of time and probably tends to a cyclic (in space) stationary state.
    The introduction of even relatively weak coalescence, seems to have a regulating
    effect that while speeds up the emergence of this phenomenon, prevents it to persist.
    \item In the case of pure spatial coalescence, with some special choices
    of initial condition and shifted coalescence kernel, the system seems to
    approach some nontrivial stationary states. Such behavior is highly vulnerable
    to addition of jumps to the system that seems to switch its asymptotic state
    to null density.
    \item For coalescence with preservation of mass the choice of cutted homogeneous
    initial condition leads to emergence of persisting irregular pattern.
    The introduction of jumps to the system seems to have a smoothing effect on those
    irregularities that makes them quickly vanishing in the flow of time.
\end{enumerate}

The introduction of both jumps and coalescence to the system seems to have
a regulating effect on each described phenomenon. One can suspect that the
observations apply to the original microscopic model as well, supposing that
the evolution of its states is well approximated by Poisson measures.

\vskip.2cm

\paragraph{\bf
Acknowledgement} Yuri Kozitsky was supported by National Science
Centre, Poland, grant 2017/25/B/ST1/00051 that is cordially
acknowledged by him.

\end{document}

%% file: fig1.tex
\begingroup
  \makeatletter
  \providecommand\color[2][]{%
    \GenericError{(gnuplot) \space\space\space\@spaces}{%
      Package color not loaded in conjunction with
      terminal option `colourtext'%
    }{See the gnuplot documentation for explanation.%
    }{Either use 'blacktext' in gnuplot or load the package
      color.sty in LaTeX.}%
    \renewcommand\color[2][]{}%
  }%
  \providecommand\includegraphics[2][]{%
    \GenericError{(gnuplot) \space\space\space\@spaces}{%
      Package graphicx or graphics not loaded%
    }{See the gnuplot documentation for explanation.%
    }{The gnuplot epslatex terminal needs graphicx.sty or graphics.sty.}%
    \renewcommand\includegraphics[2][]{}%
  }%
  \providecommand\rotatebox[2]{#2}%
  \@ifundefined{ifGPcolor}{%
    \newif\ifGPcolor
    \GPcolorfalse
  }{}%
  \@ifundefined{ifGPblacktext}{%
    \newif\ifGPblacktext
    \GPblacktexttrue
  }{}%
  \let\gplgaddtomacro\g@addto@macro
  \gdef\gplbacktext{}%
  \gdef\gplfronttext{}%
  \makeatother
  \ifGPblacktext
    \def\colorrgb#1{}%
    \def\colorgray#1{}%
  \else
    \ifGPcolor
      \def\colorrgb#1{\color[rgb]{#1}}%
      \def\colorgray#1{\color[gray]{#1}}%
      \expandafter\def\csname LTw\endcsname{\color{white}}%
      \expandafter\def\csname LTb\endcsname{\color{black}}%
      \expandafter\def\csname LTa\endcsname{\color{black}}%
      \expandafter\def\csname LT0\endcsname{\color[rgb]{1,0,0}}%
      \expandafter\def\csname LT1\endcsname{\color[rgb]{0,1,0}}%
      \expandafter\def\csname LT2\endcsname{\color[rgb]{0,0,1}}%
      \expandafter\def\csname LT3\endcsname{\color[rgb]{1,0,1}}%
      \expandafter\def\csname LT4\endcsname{\color[rgb]{0,1,1}}%
      \expandafter\def\csname LT5\endcsname{\color[rgb]{1,1,0}}%
      \expandafter\def\csname LT6\endcsname{\color[rgb]{0,0,0}}%
      \expandafter\def\csname LT7\endcsname{\color[rgb]{1,0.3,0}}%
      \expandafter\def\csname LT8\endcsname{\color[rgb]{0.5,0.5,0.5}}%
    \else
      \def\colorrgb#1{\color{black}}%
      \def\colorgray#1{\color[gray]{#1}}%
      \expandafter\def\csname LTw\endcsname{\color{white}}%
      \expandafter\def\csname LTb\endcsname{\color{black}}%
      \expandafter\def\csname LTa\endcsname{\color{black}}%
      \expandafter\def\csname LT0\endcsname{\color{black}}%
      \expandafter\def\csname LT1\endcsname{\color{black}}%
      \expandafter\def\csname LT2\endcsname{\color{black}}%
      \expandafter\def\csname LT3\endcsname{\color{black}}%
      \expandafter\def\csname LT4\endcsname{\color{black}}%
      \expandafter\def\csname LT5\endcsname{\color{black}}%
      \expandafter\def\csname LT6\endcsname{\color{black}}%
      \expandafter\def\csname LT7\endcsname{\color{black}}%
      \expandafter\def\csname LT8\endcsname{\color{black}}%
    \fi
  \fi
    \setlength{\unitlength}{0.0500bp}%
    \ifx\gptboxheight\undefined%
      \newlength{\gptboxheight}%
      \newlength{\gptboxwidth}%
      \newsavebox{\gptboxtext}%
    \fi%
    \setlength{\fboxrule}{0.5pt}%
    \setlength{\fboxsep}{1pt}%
\begin{picture}(5668.00,5668.00)%
    \gplgaddtomacro\gplbacktext{%
      \csname LTb\endcsname
      \put(594,440){\makebox(0,0)[r]{\strut{}$0$}}%
      \put(594,1394){\makebox(0,0)[r]{\strut{}$0.1$}}%
      \put(594,2347){\makebox(0,0)[r]{\strut{}$0.2$}}%
      \put(594,3301){\makebox(0,0)[r]{\strut{}$0.3$}}%
      \put(594,4255){\makebox(0,0)[r]{\strut{}$0.4$}}%
      \put(594,5209){\makebox(0,0)[r]{\strut{}$0.5$}}%
      \put(726,220){\makebox(0,0){\strut{}$-10$}}%
      \put(1862,220){\makebox(0,0){\strut{}$-5$}}%
      \put(2999,220){\makebox(0,0){\strut{}$0$}}%
      \put(4135,220){\makebox(0,0){\strut{}$5$}}%
      \put(5271,220){\makebox(0,0){\strut{}$10$}}%
    }%
    \gplgaddtomacro\gplfronttext{%
      \csname LTb\endcsname
      \put(4284,5252){\makebox(0,0)[r]{\strut{}T = 0}}%
      \csname LTb\endcsname
      \put(4284,4988){\makebox(0,0)[r]{\strut{}T = 2}}%
      \csname LTb\endcsname
      \put(4284,4724){\makebox(0,0)[r]{\strut{}T = 16}}%
      \csname LTb\endcsname
      \put(4284,4460){\makebox(0,0)[r]{\strut{}T = 80}}%
    }%
    \gplbacktext
    \put(0,0){\includegraphics{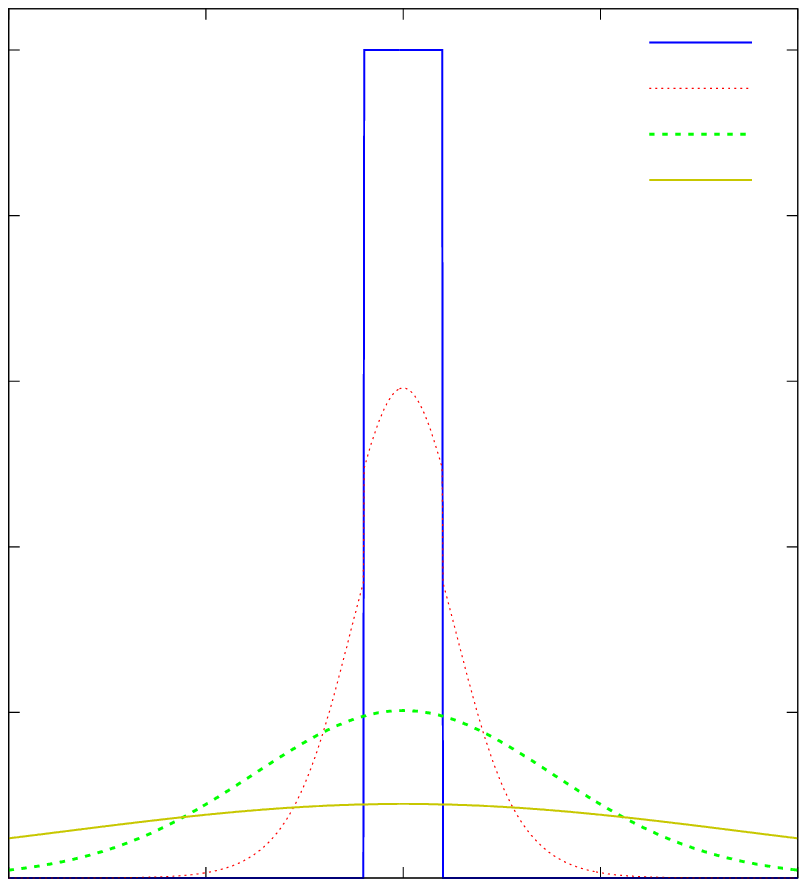}}%
    \gplfronttext
  \end{picture}%
\endgroup

%% file: fig2.tex
\begingroup
  \makeatletter
  \providecommand\color[2][]{%
    \GenericError{(gnuplot) \space\space\space\@spaces}{%
      Package color not loaded in conjunction with
      terminal option `colourtext'%
    }{See the gnuplot documentation for explanation.%
    }{Either use 'blacktext' in gnuplot or load the package
      color.sty in LaTeX.}%
    \renewcommand\color[2][]{}%
  }%
  \providecommand\includegraphics[2][]{%
    \GenericError{(gnuplot) \space\space\space\@spaces}{%
      Package graphicx or graphics not loaded%
    }{See the gnuplot documentation for explanation.%
    }{The gnuplot epslatex terminal needs graphicx.sty or graphics.sty.}%
    \renewcommand\includegraphics[2][]{}%
  }%
  \providecommand\rotatebox[2]{#2}%
  \@ifundefined{ifGPcolor}{%
    \newif\ifGPcolor
    \GPcolorfalse
  }{}%
  \@ifundefined{ifGPblacktext}{%
    \newif\ifGPblacktext
    \GPblacktexttrue
  }{}%
  \let\gplgaddtomacro\g@addto@macro
  \gdef\gplbacktext{}%
  \gdef\gplfronttext{}%
  \makeatother
  \ifGPblacktext
    \def\colorrgb#1{}%
    \def\colorgray#1{}%
  \else
    \ifGPcolor
      \def\colorrgb#1{\color[rgb]{#1}}%
      \def\colorgray#1{\color[gray]{#1}}%
      \expandafter\def\csname LTw\endcsname{\color{white}}%
      \expandafter\def\csname LTb\endcsname{\color{black}}%
      \expandafter\def\csname LTa\endcsname{\color{black}}%
      \expandafter\def\csname LT0\endcsname{\color[rgb]{1,0,0}}%
      \expandafter\def\csname LT1\endcsname{\color[rgb]{0,1,0}}%
      \expandafter\def\csname LT2\endcsname{\color[rgb]{0,0,1}}%
      \expandafter\def\csname LT3\endcsname{\color[rgb]{1,0,1}}%
      \expandafter\def\csname LT4\endcsname{\color[rgb]{0,1,1}}%
      \expandafter\def\csname LT5\endcsname{\color[rgb]{1,1,0}}%
      \expandafter\def\csname LT6\endcsname{\color[rgb]{0,0,0}}%
      \expandafter\def\csname LT7\endcsname{\color[rgb]{1,0.3,0}}%
      \expandafter\def\csname LT8\endcsname{\color[rgb]{0.5,0.5,0.5}}%
    \else
      \def\colorrgb#1{\color{black}}%
      \def\colorgray#1{\color[gray]{#1}}%
      \expandafter\def\csname LTw\endcsname{\color{white}}%
      \expandafter\def\csname LTb\endcsname{\color{black}}%
      \expandafter\def\csname LTa\endcsname{\color{black}}%
      \expandafter\def\csname LT0\endcsname{\color{black}}%
      \expandafter\def\csname LT1\endcsname{\color{black}}%
      \expandafter\def\csname LT2\endcsname{\color{black}}%
      \expandafter\def\csname LT3\endcsname{\color{black}}%
      \expandafter\def\csname LT4\endcsname{\color{black}}%
      \expandafter\def\csname LT5\endcsname{\color{black}}%
      \expandafter\def\csname LT6\endcsname{\color{black}}%
      \expandafter\def\csname LT7\endcsname{\color{black}}%
      \expandafter\def\csname LT8\endcsname{\color{black}}%
    \fi
  \fi
    \setlength{\unitlength}{0.0500bp}%
    \ifx\gptboxheight\undefined%
      \newlength{\gptboxheight}%
      \newlength{\gptboxwidth}%
      \newsavebox{\gptboxtext}%
    \fi%
    \setlength{\fboxrule}{0.5pt}%
    \setlength{\fboxsep}{1pt}%
\begin{picture}(5668.00,5668.00)%
    \gplgaddtomacro\gplbacktext{%
      \csname LTb\endcsname
      \put(726,583){\makebox(0,0)[r]{\strut{}$0.01$}}%
      \put(726,1159){\makebox(0,0)[r]{\strut{}$0.02$}}%
      \put(726,1736){\makebox(0,0)[r]{\strut{}$0.03$}}%
      \put(726,2313){\makebox(0,0)[r]{\strut{}$0.04$}}%
      \put(726,2890){\makebox(0,0)[r]{\strut{}$0.05$}}%
      \put(726,3467){\makebox(0,0)[r]{\strut{}$0.06$}}%
      \put(726,4044){\makebox(0,0)[r]{\strut{}$0.07$}}%
      \put(726,4621){\makebox(0,0)[r]{\strut{}$0.08$}}%
      \put(726,5198){\makebox(0,0)[r]{\strut{}$0.09$}}%
      \put(858,220){\makebox(0,0){\strut{}$-10$}}%
      \put(1961,220){\makebox(0,0){\strut{}$-5$}}%
      \put(3065,220){\makebox(0,0){\strut{}$0$}}%
      \put(4168,220){\makebox(0,0){\strut{}$5$}}%
      \put(5271,220){\makebox(0,0){\strut{}$10$}}%
    }%
    \gplgaddtomacro\gplfronttext{%
      \csname LTb\endcsname
      \put(4284,5252){\makebox(0,0)[r]{\strut{}$G_{1,1}$}}%
      \csname LTb\endcsname
      \put(4284,4988){\makebox(0,0)[r]{\strut{}$G_{3,1}$}}%
      \csname LTb\endcsname
      \put(4284,4724){\makebox(0,0)[r]{\strut{}$G_{1,3}$}}%
      \csname LTb\endcsname
      \put(4284,4460){\makebox(0,0)[r]{\strut{}$G_{1,1,3}$}}%
    }%
    \gplbacktext
    \put(0,0){\includegraphics{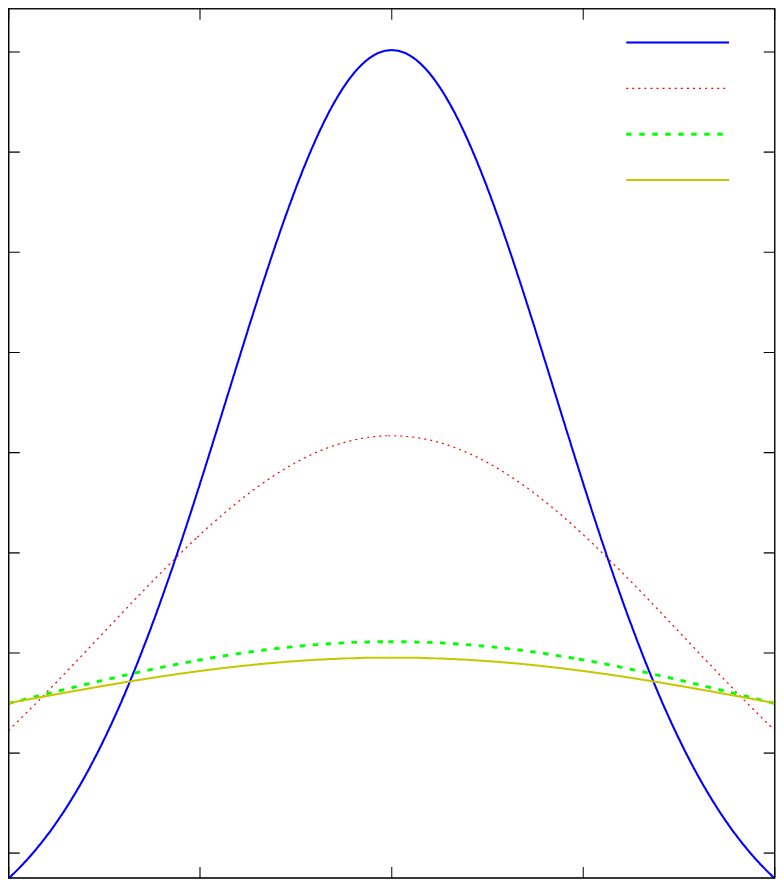}}%
    \gplfronttext
  \end{picture}%
\endgroup

%% file: fig3.tex
\begingroup
  \makeatletter
  \providecommand\color[2][]{%
    \GenericError{(gnuplot) \space\space\space\@spaces}{%
      Package color not loaded in conjunction with
      terminal option `colourtext'%
    }{See the gnuplot documentation for explanation.%
    }{Either use 'blacktext' in gnuplot or load the package
      color.sty in LaTeX.}%
    \renewcommand\color[2][]{}%
  }%
  \providecommand\includegraphics[2][]{%
    \GenericError{(gnuplot) \space\space\space\@spaces}{%
      Package graphicx or graphics not loaded%
    }{See the gnuplot documentation for explanation.%
    }{The gnuplot epslatex terminal needs graphicx.sty or graphics.sty.}%
    \renewcommand\includegraphics[2][]{}%
  }%
  \providecommand\rotatebox[2]{#2}%
  \@ifundefined{ifGPcolor}{%
    \newif\ifGPcolor
    \GPcolorfalse
  }{}%
  \@ifundefined{ifGPblacktext}{%
    \newif\ifGPblacktext
    \GPblacktexttrue
  }{}%
  \let\gplgaddtomacro\g@addto@macro
  \gdef\gplbacktext{}%
  \gdef\gplfronttext{}%
  \makeatother
  \ifGPblacktext
    \def\colorrgb#1{}%
    \def\colorgray#1{}%
  \else
    \ifGPcolor
      \def\colorrgb#1{\color[rgb]{#1}}%
      \def\colorgray#1{\color[gray]{#1}}%
      \expandafter\def\csname LTw\endcsname{\color{white}}%
      \expandafter\def\csname LTb\endcsname{\color{black}}%
      \expandafter\def\csname LTa\endcsname{\color{black}}%
      \expandafter\def\csname LT0\endcsname{\color[rgb]{1,0,0}}%
      \expandafter\def\csname LT1\endcsname{\color[rgb]{0,1,0}}%
      \expandafter\def\csname LT2\endcsname{\color[rgb]{0,0,1}}%
      \expandafter\def\csname LT3\endcsname{\color[rgb]{1,0,1}}%
      \expandafter\def\csname LT4\endcsname{\color[rgb]{0,1,1}}%
      \expandafter\def\csname LT5\endcsname{\color[rgb]{1,1,0}}%
      \expandafter\def\csname LT6\endcsname{\color[rgb]{0,0,0}}%
      \expandafter\def\csname LT7\endcsname{\color[rgb]{1,0.3,0}}%
      \expandafter\def\csname LT8\endcsname{\color[rgb]{0.5,0.5,0.5}}%
    \else
      \def\colorrgb#1{\color{black}}%
      \def\colorgray#1{\color[gray]{#1}}%
      \expandafter\def\csname LTw\endcsname{\color{white}}%
      \expandafter\def\csname LTb\endcsname{\color{black}}%
      \expandafter\def\csname LTa\endcsname{\color{black}}%
      \expandafter\def\csname LT0\endcsname{\color{black}}%
      \expandafter\def\csname LT1\endcsname{\color{black}}%
      \expandafter\def\csname LT2\endcsname{\color{black}}%
      \expandafter\def\csname LT3\endcsname{\color{black}}%
      \expandafter\def\csname LT4\endcsname{\color{black}}%
      \expandafter\def\csname LT5\endcsname{\color{black}}%
      \expandafter\def\csname LT6\endcsname{\color{black}}%
      \expandafter\def\csname LT7\endcsname{\color{black}}%
      \expandafter\def\csname LT8\endcsname{\color{black}}%
    \fi
  \fi
    \setlength{\unitlength}{0.0500bp}%
    \ifx\gptboxheight\undefined%
      \newlength{\gptboxheight}%
      \newlength{\gptboxwidth}%
      \newsavebox{\gptboxtext}%
    \fi%
    \setlength{\fboxrule}{0.5pt}%
    \setlength{\fboxsep}{1pt}%
\begin{picture}(5668.00,2834.00)%
    \gplgaddtomacro\gplbacktext{%
      \csname LTb\endcsname
      \put(594,440){\makebox(0,0)[r]{\strut{}$0$}}%
      \put(594,854){\makebox(0,0)[r]{\strut{}$0.2$}}%
      \put(594,1268){\makebox(0,0)[r]{\strut{}$0.4$}}%
      \put(594,1682){\makebox(0,0)[r]{\strut{}$0.6$}}%
      \put(594,2096){\makebox(0,0)[r]{\strut{}$0.8$}}%
      \put(594,2510){\makebox(0,0)[r]{\strut{}$1$}}%
      \put(726,220){\makebox(0,0){\strut{}$-10$}}%
      \put(1862,220){\makebox(0,0){\strut{}$-5$}}%
      \put(2999,220){\makebox(0,0){\strut{}$0$}}%
      \put(4135,220){\makebox(0,0){\strut{}$5$}}%
      \put(5271,220){\makebox(0,0){\strut{}$10$}}%
    }%
    \gplgaddtomacro\gplfronttext{%
      \csname LTb\endcsname
      \put(4284,2418){\makebox(0,0)[r]{\strut{}T = 0}}%
      \csname LTb\endcsname
      \put(4284,2154){\makebox(0,0)[r]{\strut{}T = 128}}%
      \csname LTb\endcsname
      \put(4284,1890){\makebox(0,0)[r]{\strut{}T = 640}}%
      \csname LTb\endcsname
      \put(4284,1626){\makebox(0,0)[r]{\strut{}T = 2560}}%
    }%
    \gplbacktext
    \put(0,0){\includegraphics{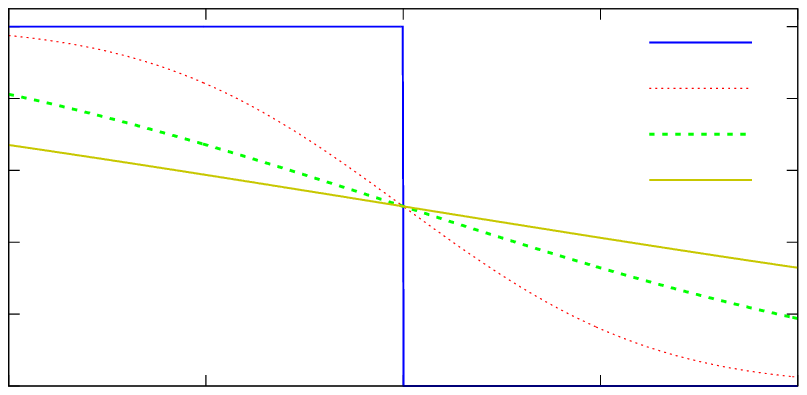}}%
    \gplfronttext
  \end{picture}%
\endgroup

%% file: fig4.tex
\begingroup
  \makeatletter
  \providecommand\color[2][]{%
    \GenericError{(gnuplot) \space\space\space\@spaces}{%
      Package color not loaded in conjunction with
      terminal option `colourtext'%
    }{See the gnuplot documentation for explanation.%
    }{Either use 'blacktext' in gnuplot or load the package
      color.sty in LaTeX.}%
    \renewcommand\color[2][]{}%
  }%
  \providecommand\includegraphics[2][]{%
    \GenericError{(gnuplot) \space\space\space\@spaces}{%
      Package graphicx or graphics not loaded%
    }{See the gnuplot documentation for explanation.%
    }{The gnuplot epslatex terminal needs graphicx.sty or graphics.sty.}%
    \renewcommand\includegraphics[2][]{}%
  }%
  \providecommand\rotatebox[2]{#2}%
  \@ifundefined{ifGPcolor}{%
    \newif\ifGPcolor
    \GPcolorfalse
  }{}%
  \@ifundefined{ifGPblacktext}{%
    \newif\ifGPblacktext
    \GPblacktexttrue
  }{}%
  \let\gplgaddtomacro\g@addto@macro
  \gdef\gplbacktext{}%
  \gdef\gplfronttext{}%
  \makeatother
  \ifGPblacktext
    \def\colorrgb#1{}%
    \def\colorgray#1{}%
  \else
    \ifGPcolor
      \def\colorrgb#1{\color[rgb]{#1}}%
      \def\colorgray#1{\color[gray]{#1}}%
      \expandafter\def\csname LTw\endcsname{\color{white}}%
      \expandafter\def\csname LTb\endcsname{\color{black}}%
      \expandafter\def\csname LTa\endcsname{\color{black}}%
      \expandafter\def\csname LT0\endcsname{\color[rgb]{1,0,0}}%
      \expandafter\def\csname LT1\endcsname{\color[rgb]{0,1,0}}%
      \expandafter\def\csname LT2\endcsname{\color[rgb]{0,0,1}}%
      \expandafter\def\csname LT3\endcsname{\color[rgb]{1,0,1}}%
      \expandafter\def\csname LT4\endcsname{\color[rgb]{0,1,1}}%
      \expandafter\def\csname LT5\endcsname{\color[rgb]{1,1,0}}%
      \expandafter\def\csname LT6\endcsname{\color[rgb]{0,0,0}}%
      \expandafter\def\csname LT7\endcsname{\color[rgb]{1,0.3,0}}%
      \expandafter\def\csname LT8\endcsname{\color[rgb]{0.5,0.5,0.5}}%
    \else
      \def\colorrgb#1{\color{black}}%
      \def\colorgray#1{\color[gray]{#1}}%
      \expandafter\def\csname LTw\endcsname{\color{white}}%
      \expandafter\def\csname LTb\endcsname{\color{black}}%
      \expandafter\def\csname LTa\endcsname{\color{black}}%
      \expandafter\def\csname LT0\endcsname{\color{black}}%
      \expandafter\def\csname LT1\endcsname{\color{black}}%
      \expandafter\def\csname LT2\endcsname{\color{black}}%
      \expandafter\def\csname LT3\endcsname{\color{black}}%
      \expandafter\def\csname LT4\endcsname{\color{black}}%
      \expandafter\def\csname LT5\endcsname{\color{black}}%
      \expandafter\def\csname LT6\endcsname{\color{black}}%
      \expandafter\def\csname LT7\endcsname{\color{black}}%
      \expandafter\def\csname LT8\endcsname{\color{black}}%
    \fi
  \fi
    \setlength{\unitlength}{0.0500bp}%
    \ifx\gptboxheight\undefined%
      \newlength{\gptboxheight}%
      \newlength{\gptboxwidth}%
      \newsavebox{\gptboxtext}%
    \fi%
    \setlength{\fboxrule}{0.5pt}%
    \setlength{\fboxsep}{1pt}%
\begin{picture}(5668.00,2834.00)%
    \gplgaddtomacro\gplbacktext{%
      \csname LTb\endcsname
      \put(594,440){\makebox(0,0)[r]{\strut{}$0$}}%
      \put(594,881){\makebox(0,0)[r]{\strut{}$0.5$}}%
      \put(594,1322){\makebox(0,0)[r]{\strut{}$1$}}%
      \put(594,1763){\makebox(0,0)[r]{\strut{}$1.5$}}%
      \put(594,2204){\makebox(0,0)[r]{\strut{}$2$}}%
      \put(726,220){\makebox(0,0){\strut{}$-40$}}%
      \put(1294,220){\makebox(0,0){\strut{}$-30$}}%
      \put(1862,220){\makebox(0,0){\strut{}$-20$}}%
      \put(2430,220){\makebox(0,0){\strut{}$-10$}}%
      \put(2999,220){\makebox(0,0){\strut{}$0$}}%
      \put(3567,220){\makebox(0,0){\strut{}$10$}}%
      \put(4135,220){\makebox(0,0){\strut{}$20$}}%
      \put(4703,220){\makebox(0,0){\strut{}$30$}}%
      \put(5271,220){\makebox(0,0){\strut{}$40$}}%
    }%
    \gplgaddtomacro\gplfronttext{%
      \csname LTb\endcsname
      \put(4284,2418){\makebox(0,0)[r]{\strut{}T = 0}}%
      \csname LTb\endcsname
      \put(4284,2154){\makebox(0,0)[r]{\strut{}T = 512}}%
      \csname LTb\endcsname
      \put(4284,1890){\makebox(0,0)[r]{\strut{}T = 1664}}%
      \csname LTb\endcsname
      \put(4284,1626){\makebox(0,0)[r]{\strut{}T = 2560}}%
    }%
    \gplbacktext
    \put(0,0){\includegraphics{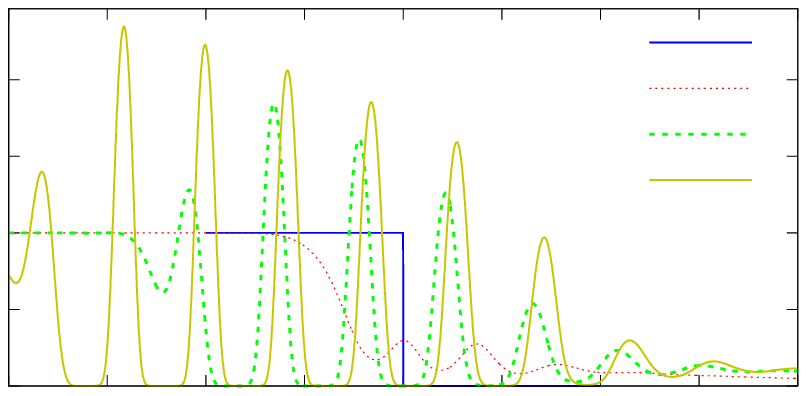}}%
    \gplfronttext
  \end{picture}%
\endgroup

%% file: fig6.tex
\begingroup
  \makeatletter
  \providecommand\color[2][]{%
    \GenericError{(gnuplot) \space\space\space\@spaces}{%
      Package color not loaded in conjunction with
      terminal option `colourtext'%
    }{See the gnuplot documentation for explanation.%
    }{Either use 'blacktext' in gnuplot or load the package
      color.sty in LaTeX.}%
    \renewcommand\color[2][]{}%
  }%
  \providecommand\includegraphics[2][]{%
    \GenericError{(gnuplot) \space\space\space\@spaces}{%
      Package graphicx or graphics not loaded%
    }{See the gnuplot documentation for explanation.%
    }{The gnuplot epslatex terminal needs graphicx.sty or graphics.sty.}%
    \renewcommand\includegraphics[2][]{}%
  }%
  \providecommand\rotatebox[2]{#2}%
  \@ifundefined{ifGPcolor}{%
    \newif\ifGPcolor
    \GPcolorfalse
  }{}%
  \@ifundefined{ifGPblacktext}{%
    \newif\ifGPblacktext
    \GPblacktexttrue
  }{}%
  \let\gplgaddtomacro\g@addto@macro
  \gdef\gplbacktext{}%
  \gdef\gplfronttext{}%
  \makeatother
  \ifGPblacktext
    \def\colorrgb#1{}%
    \def\colorgray#1{}%
  \else
    \ifGPcolor
      \def\colorrgb#1{\color[rgb]{#1}}%
      \def\colorgray#1{\color[gray]{#1}}%
      \expandafter\def\csname LTw\endcsname{\color{white}}%
      \expandafter\def\csname LTb\endcsname{\color{black}}%
      \expandafter\def\csname LTa\endcsname{\color{black}}%
      \expandafter\def\csname LT0\endcsname{\color[rgb]{1,0,0}}%
      \expandafter\def\csname LT1\endcsname{\color[rgb]{0,1,0}}%
      \expandafter\def\csname LT2\endcsname{\color[rgb]{0,0,1}}%
      \expandafter\def\csname LT3\endcsname{\color[rgb]{1,0,1}}%
      \expandafter\def\csname LT4\endcsname{\color[rgb]{0,1,1}}%
      \expandafter\def\csname LT5\endcsname{\color[rgb]{1,1,0}}%
      \expandafter\def\csname LT6\endcsname{\color[rgb]{0,0,0}}%
      \expandafter\def\csname LT7\endcsname{\color[rgb]{1,0.3,0}}%
      \expandafter\def\csname LT8\endcsname{\color[rgb]{0.5,0.5,0.5}}%
    \else
      \def\colorrgb#1{\color{black}}%
      \def\colorgray#1{\color[gray]{#1}}%
      \expandafter\def\csname LTw\endcsname{\color{white}}%
      \expandafter\def\csname LTb\endcsname{\color{black}}%
      \expandafter\def\csname LTa\endcsname{\color{black}}%
      \expandafter\def\csname LT0\endcsname{\color{black}}%
      \expandafter\def\csname LT1\endcsname{\color{black}}%
      \expandafter\def\csname LT2\endcsname{\color{black}}%
      \expandafter\def\csname LT3\endcsname{\color{black}}%
      \expandafter\def\csname LT4\endcsname{\color{black}}%
      \expandafter\def\csname LT5\endcsname{\color{black}}%
      \expandafter\def\csname LT6\endcsname{\color{black}}%
      \expandafter\def\csname LT7\endcsname{\color{black}}%
      \expandafter\def\csname LT8\endcsname{\color{black}}%
    \fi
  \fi
    \setlength{\unitlength}{0.0500bp}%
    \ifx\gptboxheight\undefined%
      \newlength{\gptboxheight}%
      \newlength{\gptboxwidth}%
      \newsavebox{\gptboxtext}%
    \fi%
    \setlength{\fboxrule}{0.5pt}%
    \setlength{\fboxsep}{1pt}%
\begin{picture}(5668.00,2834.00)%
    \gplgaddtomacro\gplbacktext{%
      \csname LTb\endcsname
      \put(594,440){\makebox(0,0)[r]{\strut{}$0$}}%
      \put(594,699){\makebox(0,0)[r]{\strut{}$0.5$}}%
      \put(594,957){\makebox(0,0)[r]{\strut{}$1$}}%
      \put(594,1216){\makebox(0,0)[r]{\strut{}$1.5$}}%
      \put(594,1475){\makebox(0,0)[r]{\strut{}$2$}}%
      \put(594,1733){\makebox(0,0)[r]{\strut{}$2.5$}}%
      \put(594,1992){\makebox(0,0)[r]{\strut{}$3$}}%
      \put(594,2251){\makebox(0,0)[r]{\strut{}$3.5$}}%
      \put(594,2510){\makebox(0,0)[r]{\strut{}$4$}}%
      \put(965,220){\makebox(0,0){\strut{}$-6$}}%
      \put(1444,220){\makebox(0,0){\strut{}$-4$}}%
      \put(1922,220){\makebox(0,0){\strut{}$-2$}}%
      \put(2400,220){\makebox(0,0){\strut{}$0$}}%
      \put(2879,220){\makebox(0,0){\strut{}$2$}}%
      \put(3357,220){\makebox(0,0){\strut{}$4$}}%
      \put(3836,220){\makebox(0,0){\strut{}$6$}}%
      \put(4314,220){\makebox(0,0){\strut{}$8$}}%
      \put(4793,220){\makebox(0,0){\strut{}$10$}}%
      \put(5271,220){\makebox(0,0){\strut{}$12$}}%
    }%
    \gplgaddtomacro\gplfronttext{%
      \csname LTb\endcsname
      \put(2373,2418){\makebox(0,0)[r]{\strut{}T = 0}}%
      \csname LTb\endcsname
      \put(2373,2154){\makebox(0,0)[r]{\strut{}T = 2}}%
      \csname LTb\endcsname
      \put(4284,2418){\makebox(0,0)[r]{\strut{}T = 320}}%
      \csname LTb\endcsname
      \put(4284,2154){\makebox(0,0)[r]{\strut{}T = 1280}}%
    }%
    \gplbacktext
    \put(0,0){\includegraphics{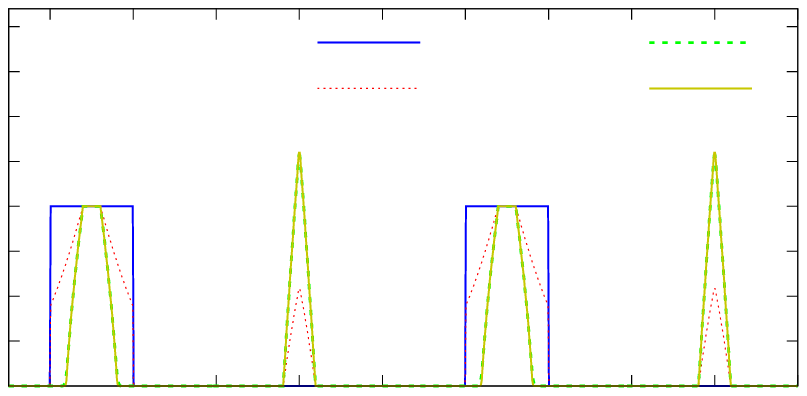}}%
    \gplfronttext
  \end{picture}%
\endgroup

%% file: fig13.tex
\begingroup
  \makeatletter
  \providecommand\color[2][]{%
    \GenericError{(gnuplot) \space\space\space\@spaces}{%
      Package color not loaded in conjunction with
      terminal option `colourtext'%
    }{See the gnuplot documentation for explanation.%
    }{Either use 'blacktext' in gnuplot or load the package
      color.sty in LaTeX.}%
    \renewcommand\color[2][]{}%
  }%
  \providecommand\includegraphics[2][]{%
    \GenericError{(gnuplot) \space\space\space\@spaces}{%
      Package graphicx or graphics not loaded%
    }{See the gnuplot documentation for explanation.%
    }{The gnuplot epslatex terminal needs graphicx.sty or graphics.sty.}%
    \renewcommand\includegraphics[2][]{}%
  }%
  \providecommand\rotatebox[2]{#2}%
  \@ifundefined{ifGPcolor}{%
    \newif\ifGPcolor
    \GPcolorfalse
  }{}%
  \@ifundefined{ifGPblacktext}{%
    \newif\ifGPblacktext
    \GPblacktexttrue
  }{}%
  \let\gplgaddtomacro\g@addto@macro
  \gdef\gplbacktext{}%
  \gdef\gplfronttext{}%
  \makeatother
  \ifGPblacktext
    \def\colorrgb#1{}%
    \def\colorgray#1{}%
  \else
    \ifGPcolor
      \def\colorrgb#1{\color[rgb]{#1}}%
      \def\colorgray#1{\color[gray]{#1}}%
      \expandafter\def\csname LTw\endcsname{\color{white}}%
      \expandafter\def\csname LTb\endcsname{\color{black}}%
      \expandafter\def\csname LTa\endcsname{\color{black}}%
      \expandafter\def\csname LT0\endcsname{\color[rgb]{1,0,0}}%
      \expandafter\def\csname LT1\endcsname{\color[rgb]{0,1,0}}%
      \expandafter\def\csname LT2\endcsname{\color[rgb]{0,0,1}}%
      \expandafter\def\csname LT3\endcsname{\color[rgb]{1,0,1}}%
      \expandafter\def\csname LT4\endcsname{\color[rgb]{0,1,1}}%
      \expandafter\def\csname LT5\endcsname{\color[rgb]{1,1,0}}%
      \expandafter\def\csname LT6\endcsname{\color[rgb]{0,0,0}}%
      \expandafter\def\csname LT7\endcsname{\color[rgb]{1,0.3,0}}%
      \expandafter\def\csname LT8\endcsname{\color[rgb]{0.5,0.5,0.5}}%
    \else
      \def\colorrgb#1{\color{black}}%
      \def\colorgray#1{\color[gray]{#1}}%
      \expandafter\def\csname LTw\endcsname{\color{white}}%
      \expandafter\def\csname LTb\endcsname{\color{black}}%
      \expandafter\def\csname LTa\endcsname{\color{black}}%
      \expandafter\def\csname LT0\endcsname{\color{black}}%
      \expandafter\def\csname LT1\endcsname{\color{black}}%
      \expandafter\def\csname LT2\endcsname{\color{black}}%
      \expandafter\def\csname LT3\endcsname{\color{black}}%
      \expandafter\def\csname LT4\endcsname{\color{black}}%
      \expandafter\def\csname LT5\endcsname{\color{black}}%
      \expandafter\def\csname LT6\endcsname{\color{black}}%
      \expandafter\def\csname LT7\endcsname{\color{black}}%
      \expandafter\def\csname LT8\endcsname{\color{black}}%
    \fi
  \fi
    \setlength{\unitlength}{0.0500bp}%
    \ifx\gptboxheight\undefined%
      \newlength{\gptboxheight}%
      \newlength{\gptboxwidth}%
      \newsavebox{\gptboxtext}%
    \fi%
    \setlength{\fboxrule}{0.5pt}%
    \setlength{\fboxsep}{1pt}%
\begin{picture}(5668.00,5668.00)%
    \gplgaddtomacro\gplbacktext{%
      \csname LTb\endcsname
      \put(594,440){\makebox(0,0)[r]{\strut{}$0$}}%
      \put(594,1235){\makebox(0,0)[r]{\strut{}$0.2$}}%
      \put(594,2030){\makebox(0,0)[r]{\strut{}$0.4$}}%
      \put(594,2824){\makebox(0,0)[r]{\strut{}$0.6$}}%
      \put(594,3619){\makebox(0,0)[r]{\strut{}$0.8$}}%
      \put(594,4414){\makebox(0,0)[r]{\strut{}$1$}}%
      \put(594,5209){\makebox(0,0)[r]{\strut{}$1.2$}}%
      \put(726,220){\makebox(0,0){\strut{}$-1$}}%
      \put(1635,220){\makebox(0,0){\strut{}$0$}}%
      \put(2544,220){\makebox(0,0){\strut{}$1$}}%
      \put(3453,220){\makebox(0,0){\strut{}$2$}}%
      \put(4362,220){\makebox(0,0){\strut{}$3$}}%
      \put(5271,220){\makebox(0,0){\strut{}$4$}}%
    }%
    \gplgaddtomacro\gplfronttext{%
      \csname LTb\endcsname
      \put(2373,5252){\makebox(0,0)[r]{\strut{}T = 0}}%
      \csname LTb\endcsname
      \put(2373,4988){\makebox(0,0)[r]{\strut{}T = 64}}%
      \csname LTb\endcsname
      \put(4284,5252){\makebox(0,0)[r]{\strut{}T = 192}}%
      \csname LTb\endcsname
      \put(4284,4988){\makebox(0,0)[r]{\strut{}T = 1280}}%
    }%
    \gplbacktext
    \put(0,0){\includegraphics{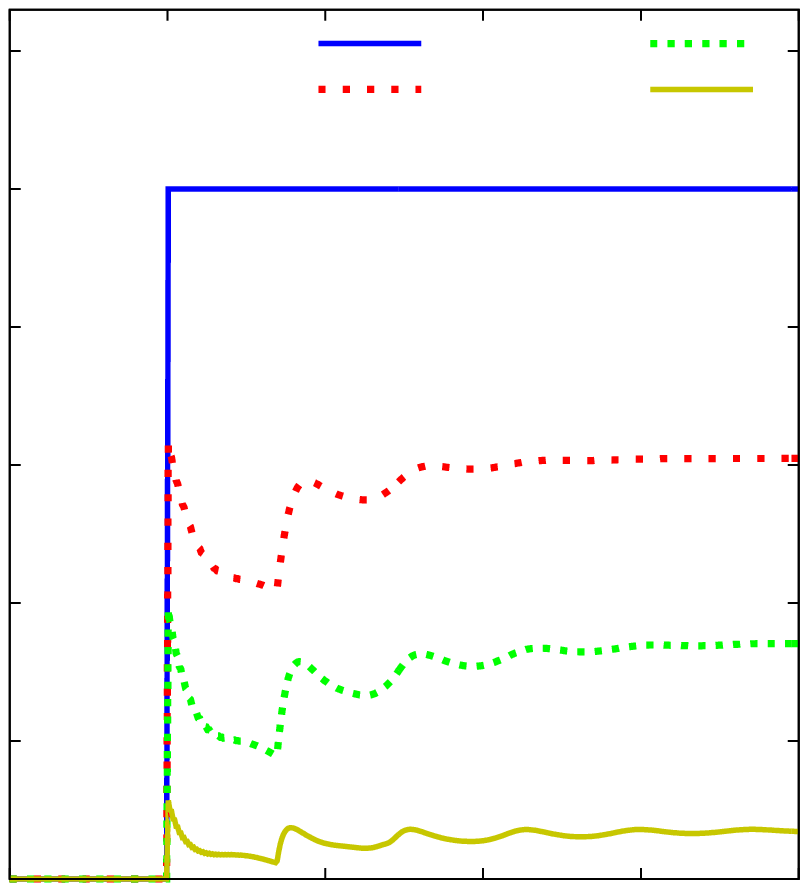}}%
    \gplfronttext
  \end{picture}%
\endgroup

%% file: fig7.tex
\begingroup
  \makeatletter
  \providecommand\color[2][]{%
    \GenericError{(gnuplot) \space\space\space\@spaces}{%
      Package color not loaded in conjunction with
      terminal option `colourtext'%
    }{See the gnuplot documentation for explanation.%
    }{Either use 'blacktext' in gnuplot or load the package
      color.sty in LaTeX.}%
    \renewcommand\color[2][]{}%
  }%
  \providecommand\includegraphics[2][]{%
    \GenericError{(gnuplot) \space\space\space\@spaces}{%
      Package graphicx or graphics not loaded%
    }{See the gnuplot documentation for explanation.%
    }{The gnuplot epslatex terminal needs graphicx.sty or graphics.sty.}%
    \renewcommand\includegraphics[2][]{}%
  }%
  \providecommand\rotatebox[2]{#2}%
  \@ifundefined{ifGPcolor}{%
    \newif\ifGPcolor
    \GPcolorfalse
  }{}%
  \@ifundefined{ifGPblacktext}{%
    \newif\ifGPblacktext
    \GPblacktexttrue
  }{}%
  \let\gplgaddtomacro\g@addto@macro
  \gdef\gplbacktext{}%
  \gdef\gplfronttext{}%
  \makeatother
  \ifGPblacktext
    \def\colorrgb#1{}%
    \def\colorgray#1{}%
  \else
    \ifGPcolor
      \def\colorrgb#1{\color[rgb]{#1}}%
      \def\colorgray#1{\color[gray]{#1}}%
      \expandafter\def\csname LTw\endcsname{\color{white}}%
      \expandafter\def\csname LTb\endcsname{\color{black}}%
      \expandafter\def\csname LTa\endcsname{\color{black}}%
      \expandafter\def\csname LT0\endcsname{\color[rgb]{1,0,0}}%
      \expandafter\def\csname LT1\endcsname{\color[rgb]{0,1,0}}%
      \expandafter\def\csname LT2\endcsname{\color[rgb]{0,0,1}}%
      \expandafter\def\csname LT3\endcsname{\color[rgb]{1,0,1}}%
      \expandafter\def\csname LT4\endcsname{\color[rgb]{0,1,1}}%
      \expandafter\def\csname LT5\endcsname{\color[rgb]{1,1,0}}%
      \expandafter\def\csname LT6\endcsname{\color[rgb]{0,0,0}}%
      \expandafter\def\csname LT7\endcsname{\color[rgb]{1,0.3,0}}%
      \expandafter\def\csname LT8\endcsname{\color[rgb]{0.5,0.5,0.5}}%
    \else
      \def\colorrgb#1{\color{black}}%
      \def\colorgray#1{\color[gray]{#1}}%
      \expandafter\def\csname LTw\endcsname{\color{white}}%
      \expandafter\def\csname LTb\endcsname{\color{black}}%
      \expandafter\def\csname LTa\endcsname{\color{black}}%
      \expandafter\def\csname LT0\endcsname{\color{black}}%
      \expandafter\def\csname LT1\endcsname{\color{black}}%
      \expandafter\def\csname LT2\endcsname{\color{black}}%
      \expandafter\def\csname LT3\endcsname{\color{black}}%
      \expandafter\def\csname LT4\endcsname{\color{black}}%
      \expandafter\def\csname LT5\endcsname{\color{black}}%
      \expandafter\def\csname LT6\endcsname{\color{black}}%
      \expandafter\def\csname LT7\endcsname{\color{black}}%
      \expandafter\def\csname LT8\endcsname{\color{black}}%
    \fi
  \fi
    \setlength{\unitlength}{0.0500bp}%
    \ifx\gptboxheight\undefined%
      \newlength{\gptboxheight}%
      \newlength{\gptboxwidth}%
      \newsavebox{\gptboxtext}%
    \fi%
    \setlength{\fboxrule}{0.5pt}%
    \setlength{\fboxsep}{1pt}%
\begin{picture}(5668.00,5668.00)%
    \gplgaddtomacro\gplbacktext{%
      \csname LTb\endcsname
      \put(594,440){\makebox(0,0)[r]{\strut{}$0$}}%
      \put(594,1235){\makebox(0,0)[r]{\strut{}$0.2$}}%
      \put(594,2030){\makebox(0,0)[r]{\strut{}$0.4$}}%
      \put(594,2824){\makebox(0,0)[r]{\strut{}$0.6$}}%
      \put(594,3619){\makebox(0,0)[r]{\strut{}$0.8$}}%
      \put(594,4414){\makebox(0,0)[r]{\strut{}$1$}}%
      \put(594,5209){\makebox(0,0)[r]{\strut{}$1.2$}}%
      \put(726,220){\makebox(0,0){\strut{}$-1$}}%
      \put(1635,220){\makebox(0,0){\strut{}$0$}}%
      \put(2544,220){\makebox(0,0){\strut{}$1$}}%
      \put(3453,220){\makebox(0,0){\strut{}$2$}}%
      \put(4362,220){\makebox(0,0){\strut{}$3$}}%
      \put(5271,220){\makebox(0,0){\strut{}$4$}}%
    }%
    \gplgaddtomacro\gplfronttext{%
      \csname LTb\endcsname
      \put(2373,5252){\makebox(0,0)[r]{\strut{}T = 0}}%
      \csname LTb\endcsname
      \put(2373,4988){\makebox(0,0)[r]{\strut{}T = 64}}%
      \csname LTb\endcsname
      \put(4284,5252){\makebox(0,0)[r]{\strut{}T = 192}}%
      \csname LTb\endcsname
      \put(4284,4988){\makebox(0,0)[r]{\strut{}T = 1280}}%
    }%
    \gplbacktext
    \put(0,0){\includegraphics{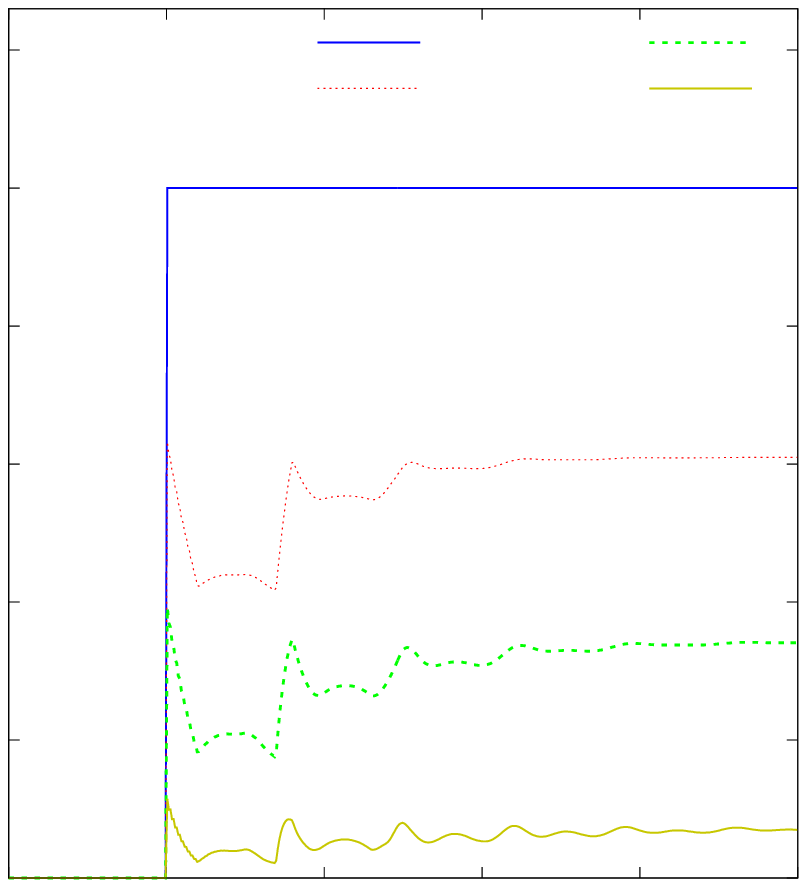}}%
    \gplfronttext
  \end{picture}%
\endgroup

%% file: fig14.tex
\begingroup
  \makeatletter
  \providecommand\color[2][]{%
    \GenericError{(gnuplot) \space\space\space\@spaces}{%
      Package color not loaded in conjunction with
      terminal option `colourtext'%
    }{See the gnuplot documentation for explanation.%
    }{Either use 'blacktext' in gnuplot or load the package
      color.sty in LaTeX.}%
    \renewcommand\color[2][]{}%
  }%
  \providecommand\includegraphics[2][]{%
    \GenericError{(gnuplot) \space\space\space\@spaces}{%
      Package graphicx or graphics not loaded%
    }{See the gnuplot documentation for explanation.%
    }{The gnuplot epslatex terminal needs graphicx.sty or graphics.sty.}%
    \renewcommand\includegraphics[2][]{}%
  }%
  \providecommand\rotatebox[2]{#2}%
  \@ifundefined{ifGPcolor}{%
    \newif\ifGPcolor
    \GPcolorfalse
  }{}%
  \@ifundefined{ifGPblacktext}{%
    \newif\ifGPblacktext
    \GPblacktexttrue
  }{}%
  \let\gplgaddtomacro\g@addto@macro
  \gdef\gplbacktext{}%
  \gdef\gplfronttext{}%
  \makeatother
  \ifGPblacktext
    \def\colorrgb#1{}%
    \def\colorgray#1{}%
  \else
    \ifGPcolor
      \def\colorrgb#1{\color[rgb]{#1}}%
      \def\colorgray#1{\color[gray]{#1}}%
      \expandafter\def\csname LTw\endcsname{\color{white}}%
      \expandafter\def\csname LTb\endcsname{\color{black}}%
      \expandafter\def\csname LTa\endcsname{\color{black}}%
      \expandafter\def\csname LT0\endcsname{\color[rgb]{1,0,0}}%
      \expandafter\def\csname LT1\endcsname{\color[rgb]{0,1,0}}%
      \expandafter\def\csname LT2\endcsname{\color[rgb]{0,0,1}}%
      \expandafter\def\csname LT3\endcsname{\color[rgb]{1,0,1}}%
      \expandafter\def\csname LT4\endcsname{\color[rgb]{0,1,1}}%
      \expandafter\def\csname LT5\endcsname{\color[rgb]{1,1,0}}%
      \expandafter\def\csname LT6\endcsname{\color[rgb]{0,0,0}}%
      \expandafter\def\csname LT7\endcsname{\color[rgb]{1,0.3,0}}%
      \expandafter\def\csname LT8\endcsname{\color[rgb]{0.5,0.5,0.5}}%
    \else
      \def\colorrgb#1{\color{black}}%
      \def\colorgray#1{\color[gray]{#1}}%
      \expandafter\def\csname LTw\endcsname{\color{white}}%
      \expandafter\def\csname LTb\endcsname{\color{black}}%
      \expandafter\def\csname LTa\endcsname{\color{black}}%
      \expandafter\def\csname LT0\endcsname{\color{black}}%
      \expandafter\def\csname LT1\endcsname{\color{black}}%
      \expandafter\def\csname LT2\endcsname{\color{black}}%
      \expandafter\def\csname LT3\endcsname{\color{black}}%
      \expandafter\def\csname LT4\endcsname{\color{black}}%
      \expandafter\def\csname LT5\endcsname{\color{black}}%
      \expandafter\def\csname LT6\endcsname{\color{black}}%
      \expandafter\def\csname LT7\endcsname{\color{black}}%
      \expandafter\def\csname LT8\endcsname{\color{black}}%
    \fi
  \fi
    \setlength{\unitlength}{0.0500bp}%
    \ifx\gptboxheight\undefined%
      \newlength{\gptboxheight}%
      \newlength{\gptboxwidth}%
      \newsavebox{\gptboxtext}%
    \fi%
    \setlength{\fboxrule}{0.5pt}%
    \setlength{\fboxsep}{1pt}%
\begin{picture}(7936.00,3400.00)%
    \gplgaddtomacro\gplbacktext{%
      \csname LTb\endcsname
      \put(594,440){\makebox(0,0)[r]{\strut{}$0$}}%
      \put(594,841){\makebox(0,0)[r]{\strut{}$0.2$}}%
      \put(594,1243){\makebox(0,0)[r]{\strut{}$0.4$}}%
      \put(594,1644){\makebox(0,0)[r]{\strut{}$0.6$}}%
      \put(594,2045){\makebox(0,0)[r]{\strut{}$0.8$}}%
      \put(594,2447){\makebox(0,0)[r]{\strut{}$1$}}%
      \put(594,2848){\makebox(0,0)[r]{\strut{}$1.2$}}%
      \put(726,220){\makebox(0,0){\strut{}$-70$}}%
      \put(1578,220){\makebox(0,0){\strut{}$-60$}}%
      \put(2429,220){\makebox(0,0){\strut{}$-50$}}%
      \put(3281,220){\makebox(0,0){\strut{}$-40$}}%
      \put(4133,220){\makebox(0,0){\strut{}$-30$}}%
      \put(4984,220){\makebox(0,0){\strut{}$-20$}}%
      \put(5836,220){\makebox(0,0){\strut{}$-10$}}%
      \put(6687,220){\makebox(0,0){\strut{}$0$}}%
      \put(7539,220){\makebox(0,0){\strut{}$10$}}%
    }%
    \gplgaddtomacro\gplfronttext{%
      \csname LTb\endcsname
      \put(2994,2984){\makebox(0,0)[r]{\strut{}T = 0}}%
      \csname LTb\endcsname
      \put(2994,2720){\makebox(0,0)[r]{\strut{}T = 8}}%
      \csname LTb\endcsname
      \put(4773,2984){\makebox(0,0)[r]{\strut{}T = 32}}%
      \csname LTb\endcsname
      \put(4773,2720){\makebox(0,0)[r]{\strut{}T = 192}}%
      \csname LTb\endcsname
      \put(6552,2984){\makebox(0,0)[r]{\strut{}T = 256}}%
      \csname LTb\endcsname
      \put(6552,2720){\makebox(0,0)[r]{\strut{}T = 320}}%
    }%
    \gplbacktext
    \put(0,0){\includegraphics{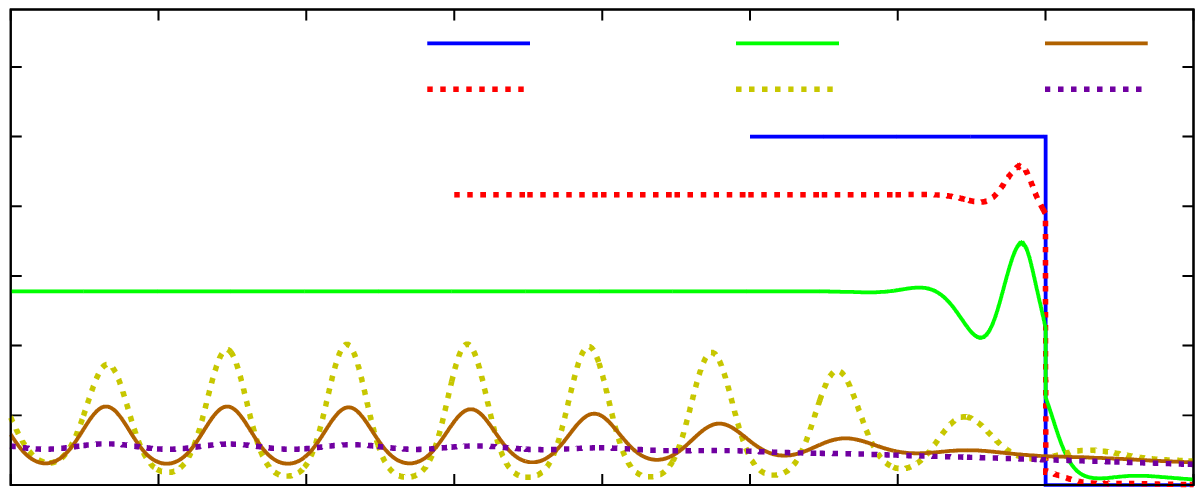}}%
    \gplfronttext
  \end{picture}%
\endgroup

%% file: fig10.tex
\begingroup
  \makeatletter
  \providecommand\color[2][]{%
    \GenericError{(gnuplot) \space\space\space\@spaces}{%
      Package color not loaded in conjunction with
      terminal option `colourtext'%
    }{See the gnuplot documentation for explanation.%
    }{Either use 'blacktext' in gnuplot or load the package
      color.sty in LaTeX.}%
    \renewcommand\color[2][]{}%
  }%
  \providecommand\includegraphics[2][]{%
    \GenericError{(gnuplot) \space\space\space\@spaces}{%
      Package graphicx or graphics not loaded%
    }{See the gnuplot documentation for explanation.%
    }{The gnuplot epslatex terminal needs graphicx.sty or graphics.sty.}%
    \renewcommand\includegraphics[2][]{}%
  }%
  \providecommand\rotatebox[2]{#2}%
  \@ifundefined{ifGPcolor}{%
    \newif\ifGPcolor
    \GPcolorfalse
  }{}%
  \@ifundefined{ifGPblacktext}{%
    \newif\ifGPblacktext
    \GPblacktexttrue
  }{}%
  \let\gplgaddtomacro\g@addto@macro
  \gdef\gplbacktext{}%
  \gdef\gplfronttext{}%
  \makeatother
  \ifGPblacktext
    \def\colorrgb#1{}%
    \def\colorgray#1{}%
  \else
    \ifGPcolor
      \def\colorrgb#1{\color[rgb]{#1}}%
      \def\colorgray#1{\color[gray]{#1}}%
      \expandafter\def\csname LTw\endcsname{\color{white}}%
      \expandafter\def\csname LTb\endcsname{\color{black}}%
      \expandafter\def\csname LTa\endcsname{\color{black}}%
      \expandafter\def\csname LT0\endcsname{\color[rgb]{1,0,0}}%
      \expandafter\def\csname LT1\endcsname{\color[rgb]{0,1,0}}%
      \expandafter\def\csname LT2\endcsname{\color[rgb]{0,0,1}}%
      \expandafter\def\csname LT3\endcsname{\color[rgb]{1,0,1}}%
      \expandafter\def\csname LT4\endcsname{\color[rgb]{0,1,1}}%
      \expandafter\def\csname LT5\endcsname{\color[rgb]{1,1,0}}%
      \expandafter\def\csname LT6\endcsname{\color[rgb]{0,0,0}}%
      \expandafter\def\csname LT7\endcsname{\color[rgb]{1,0.3,0}}%
      \expandafter\def\csname LT8\endcsname{\color[rgb]{0.5,0.5,0.5}}%
    \else
      \def\colorrgb#1{\color{black}}%
      \def\colorgray#1{\color[gray]{#1}}%
      \expandafter\def\csname LTw\endcsname{\color{white}}%
      \expandafter\def\csname LTb\endcsname{\color{black}}%
      \expandafter\def\csname LTa\endcsname{\color{black}}%
      \expandafter\def\csname LT0\endcsname{\color{black}}%
      \expandafter\def\csname LT1\endcsname{\color{black}}%
      \expandafter\def\csname LT2\endcsname{\color{black}}%
      \expandafter\def\csname LT3\endcsname{\color{black}}%
      \expandafter\def\csname LT4\endcsname{\color{black}}%
      \expandafter\def\csname LT5\endcsname{\color{black}}%
      \expandafter\def\csname LT6\endcsname{\color{black}}%
      \expandafter\def\csname LT7\endcsname{\color{black}}%
      \expandafter\def\csname LT8\endcsname{\color{black}}%
    \fi
  \fi
    \setlength{\unitlength}{0.0500bp}%
    \ifx\gptboxheight\undefined%
      \newlength{\gptboxheight}%
      \newlength{\gptboxwidth}%
      \newsavebox{\gptboxtext}%
    \fi%
    \setlength{\fboxrule}{0.5pt}%
    \setlength{\fboxsep}{1pt}%
\begin{picture}(7936.00,3400.00)%
    \gplgaddtomacro\gplbacktext{%
      \csname LTb\endcsname
      \put(594,440){\makebox(0,0)[r]{\strut{}$0$}}%
      \put(594,766){\makebox(0,0)[r]{\strut{}$0.5$}}%
      \put(594,1092){\makebox(0,0)[r]{\strut{}$1$}}%
      \put(594,1418){\makebox(0,0)[r]{\strut{}$1.5$}}%
      \put(594,1744){\makebox(0,0)[r]{\strut{}$2$}}%
      \put(594,2070){\makebox(0,0)[r]{\strut{}$2.5$}}%
      \put(594,2396){\makebox(0,0)[r]{\strut{}$3$}}%
      \put(594,2723){\makebox(0,0)[r]{\strut{}$3.5$}}%
      \put(594,3049){\makebox(0,0)[r]{\strut{}$4$}}%
      \put(1085,220){\makebox(0,0){\strut{}$-6$}}%
      \put(1802,220){\makebox(0,0){\strut{}$-4$}}%
      \put(2519,220){\makebox(0,0){\strut{}$-2$}}%
      \put(3236,220){\makebox(0,0){\strut{}$0$}}%
      \put(3953,220){\makebox(0,0){\strut{}$2$}}%
      \put(4670,220){\makebox(0,0){\strut{}$4$}}%
      \put(5388,220){\makebox(0,0){\strut{}$6$}}%
      \put(6105,220){\makebox(0,0){\strut{}$8$}}%
      \put(6822,220){\makebox(0,0){\strut{}$10$}}%
      \put(7539,220){\makebox(0,0){\strut{}$12$}}%
    }%
    \gplgaddtomacro\gplfronttext{%
      \csname LTb\endcsname
      \put(4905,2984){\makebox(0,0)[r]{\strut{}T = 0}}%
      \csname LTb\endcsname
      \put(4905,2720){\makebox(0,0)[r]{\strut{}T = 2}}%
      \csname LTb\endcsname
      \put(6552,2984){\makebox(0,0)[r]{\strut{}T = 10}}%
      \csname LTb\endcsname
      \put(6552,2720){\makebox(0,0)[r]{\strut{}T = 30}}%
    }%
    \gplbacktext
    \put(0,0){\includegraphics{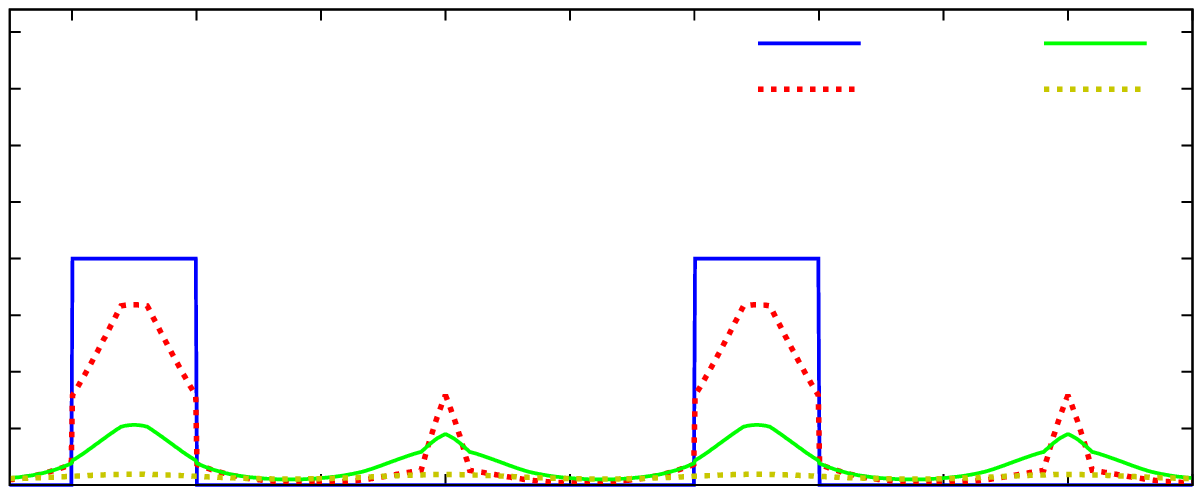}}%
    \gplfronttext
  \end{picture}%
\endgroup

%% file: fig12.tex
\begingroup
  \makeatletter
  \providecommand\color[2][]{%
    \GenericError{(gnuplot) \space\space\space\@spaces}{%
      Package color not loaded in conjunction with
      terminal option `colourtext'%
    }{See the gnuplot documentation for explanation.%
    }{Either use 'blacktext' in gnuplot or load the package
      color.sty in LaTeX.}%
    \renewcommand\color[2][]{}%
  }%
  \providecommand\includegraphics[2][]{%
    \GenericError{(gnuplot) \space\space\space\@spaces}{%
      Package graphicx or graphics not loaded%
    }{See the gnuplot documentation for explanation.%
    }{The gnuplot epslatex terminal needs graphicx.sty or graphics.sty.}%
    \renewcommand\includegraphics[2][]{}%
  }%
  \providecommand\rotatebox[2]{#2}%
  \@ifundefined{ifGPcolor}{%
    \newif\ifGPcolor
    \GPcolorfalse
  }{}%
  \@ifundefined{ifGPblacktext}{%
    \newif\ifGPblacktext
    \GPblacktexttrue
  }{}%
  \let\gplgaddtomacro\g@addto@macro
  \gdef\gplbacktext{}%
  \gdef\gplfronttext{}%
  \makeatother
  \ifGPblacktext
    \def\colorrgb#1{}%
    \def\colorgray#1{}%
  \else
    \ifGPcolor
      \def\colorrgb#1{\color[rgb]{#1}}%
      \def\colorgray#1{\color[gray]{#1}}%
      \expandafter\def\csname LTw\endcsname{\color{white}}%
      \expandafter\def\csname LTb\endcsname{\color{black}}%
      \expandafter\def\csname LTa\endcsname{\color{black}}%
      \expandafter\def\csname LT0\endcsname{\color[rgb]{1,0,0}}%
      \expandafter\def\csname LT1\endcsname{\color[rgb]{0,1,0}}%
      \expandafter\def\csname LT2\endcsname{\color[rgb]{0,0,1}}%
      \expandafter\def\csname LT3\endcsname{\color[rgb]{1,0,1}}%
      \expandafter\def\csname LT4\endcsname{\color[rgb]{0,1,1}}%
      \expandafter\def\csname LT5\endcsname{\color[rgb]{1,1,0}}%
      \expandafter\def\csname LT6\endcsname{\color[rgb]{0,0,0}}%
      \expandafter\def\csname LT7\endcsname{\color[rgb]{1,0.3,0}}%
      \expandafter\def\csname LT8\endcsname{\color[rgb]{0.5,0.5,0.5}}%
    \else
      \def\colorrgb#1{\color{black}}%
      \def\colorgray#1{\color[gray]{#1}}%
      \expandafter\def\csname LTw\endcsname{\color{white}}%
      \expandafter\def\csname LTb\endcsname{\color{black}}%
      \expandafter\def\csname LTa\endcsname{\color{black}}%
      \expandafter\def\csname LT0\endcsname{\color{black}}%
      \expandafter\def\csname LT1\endcsname{\color{black}}%
      \expandafter\def\csname LT2\endcsname{\color{black}}%
      \expandafter\def\csname LT3\endcsname{\color{black}}%
      \expandafter\def\csname LT4\endcsname{\color{black}}%
      \expandafter\def\csname LT5\endcsname{\color{black}}%
      \expandafter\def\csname LT6\endcsname{\color{black}}%
      \expandafter\def\csname LT7\endcsname{\color{black}}%
      \expandafter\def\csname LT8\endcsname{\color{black}}%
    \fi
  \fi
    \setlength{\unitlength}{0.0500bp}%
    \ifx\gptboxheight\undefined%
      \newlength{\gptboxheight}%
      \newlength{\gptboxwidth}%
      \newsavebox{\gptboxtext}%
    \fi%
    \setlength{\fboxrule}{0.5pt}%
    \setlength{\fboxsep}{1pt}%
\begin{picture}(7936.00,3400.00)%
    \gplgaddtomacro\gplbacktext{%
      \csname LTb\endcsname
      \put(594,440){\makebox(0,0)[r]{\strut{}$0$}}%
      \put(594,841){\makebox(0,0)[r]{\strut{}$0.2$}}%
      \put(594,1243){\makebox(0,0)[r]{\strut{}$0.4$}}%
      \put(594,1644){\makebox(0,0)[r]{\strut{}$0.6$}}%
      \put(594,2045){\makebox(0,0)[r]{\strut{}$0.8$}}%
      \put(594,2447){\makebox(0,0)[r]{\strut{}$1$}}%
      \put(594,2848){\makebox(0,0)[r]{\strut{}$1.2$}}%
      \put(726,220){\makebox(0,0){\strut{}$-1$}}%
      \put(2089,220){\makebox(0,0){\strut{}$0$}}%
      \put(3451,220){\makebox(0,0){\strut{}$1$}}%
      \put(4814,220){\makebox(0,0){\strut{}$2$}}%
      \put(6176,220){\makebox(0,0){\strut{}$3$}}%
      \put(7539,220){\makebox(0,0){\strut{}$4$}}%
    }%
    \gplgaddtomacro\gplfronttext{%
      \csname LTb\endcsname
      \put(4641,2984){\makebox(0,0)[r]{\strut{}T = 0}}%
      \csname LTb\endcsname
      \put(4641,2720){\makebox(0,0)[r]{\strut{}T = 64}}%
      \csname LTb\endcsname
      \put(6552,2984){\makebox(0,0)[r]{\strut{}T = 192}}%
      \csname LTb\endcsname
      \put(6552,2720){\makebox(0,0)[r]{\strut{}T = 1280}}%
    }%
    \gplbacktext
    \put(0,0){\includegraphics{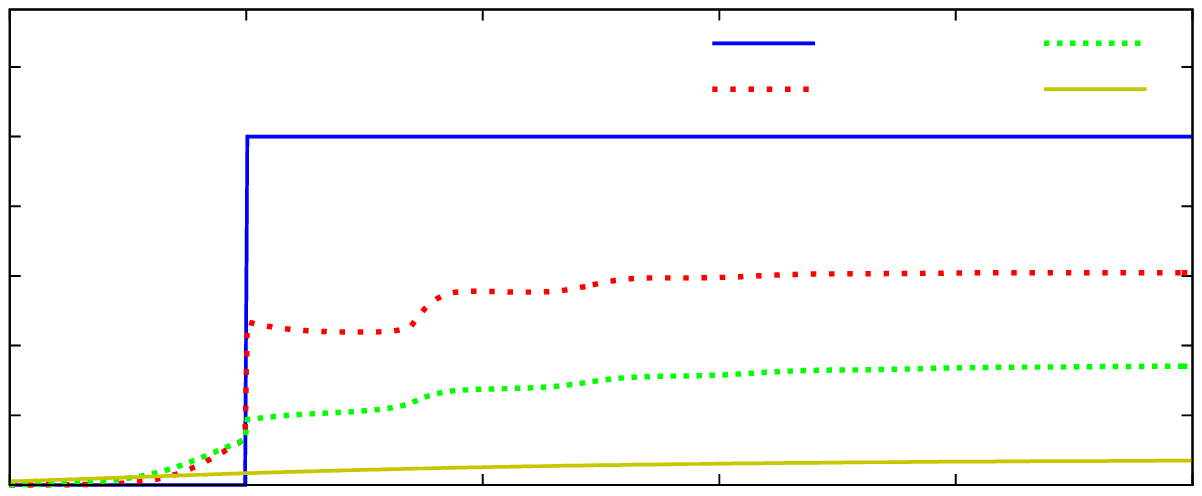}}%
    \gplfronttext
  \end{picture}%
\endgroup

%% file: jumpsAndCoalescenceFFF.bbl
\begin{thebibliography}{ll}
\bibitem{Adams} Adams, T. P., Holland, E. P.,  Law, R., Plank  M. J.,
Raghib, M.: On the growth of locally interacting plants:
differential equations for the dynamics of spatial moments, Ecology
94, 2732--2743 (2013)

\bibitem{Arrat} Arratia, R. A.: Coalescing Brownian Motion on the
Line. ProQuest LLC, Ann Arbor, MI, 1979. Thesis (Ph.D) - The
University of Wisconsin - Madison
\bibitem{Ba} Banasiak, J., Lachowicz, M.: Methods of Small Parameter
in Mathematical Biology. Birh\"auser/Springer, Cham (2014)
\bibitem{Kawasaki2} Bara\'nska, J., Kozitsky, Yu.: The global
evolution of states of a continuum Kawasaki model with repulsion.
IMA J. Appl. Math. 83, 412--435 (2018)
\bibitem{Beres} Berestycki, N., Garban, Ch., Sen, A.: Coalescing
Brownian flows: a new approach. Ann. Probab. 43, 3177--3215 (2015)

\bibitem{Bogol} Bogoliubov, N. N.:  Problems of a Dynamical Theory in Statistical Physics. In: Studies in
Statistical Mechanics (ed. J. de Boer $\&$ G. E. Uhenbeck).
North-Holland (1962)

\bibitem{BP1}  Bolker, B.M,  Pacala, S. W.:
Using moment equations to understand stochastically driven spatial
pattern formation in ecological systems, Theoret. Population Biol.
52, 179--197  (1997)

\bibitem{Delius} Capit\'an, J. A., Delius, G. W.:  
Scale-invariant model of marine population dynamics, 
Phys. Rev. E 81, 061901 (2010)

\bibitem{Cer} Cercignani, C.: On the Boltzmann equation for rigid spheres, Transport Th. and Stat.
Phys. 2,  211--225 (1972)

\bibitem{Boltz} Cercignani C.: The Boltzmann equation. In: The
Boltzmann Equation and Its Applications. Applied Mathematical
Sciences, vol 67. Springer, New York, NY (1988)

\bibitem{Dawson} Dawson, D. A.: Measure-Valued Markov Processes. {\'E}cole d'{\'E}t{\'e} de
Probabilit{\'e}s de Saint-Flour XXI--1991, 1--260, Lecture Notes in
Math., 1541, Springer, Berlin (1993)
\bibitem{Reb} Dorlas, T. C., Rebenko, A. L., Savoie, B.: Correlation of clusters: Partially truncated correlation functions and their
decay, J. Math, Phys. 61, 033303 (2020)




\bibitem{FKKK} Finkelshtein, D., Kondratiev, Yu., Kozitsky, Yu.,
Kutoviy, O.: The statistical dynamics of a spatial logistic model
and the related kinetic equation. Math. Models Methods Appl. Sci.
25, 343--370 (2015)
\bibitem{Kovn} Konarovskii, V. V.: On an infinite system of
diffusing particles with coalescing. Teor. Veroyatn. Primen. 55,
157--167 (2010)
\bibitem{KovnR} Konarovskii, V. V., von Renesse, M.: Modified
massive Arratia flow and Wasserstein diffusion. Comm. Pure Appl.
Math. 72, 764--800 (2019)
\bibitem{SpatialEcologicalModel} Kondratiev, Yu., Kozitsky, Yu.: The
evolution of states in a spatial pupulation model. J. Dynam.
Differential Equations 30, 135--173 (2018)

\bibitem{KoP} Kozitsky, Yu., Pilorz, K.: Random jumps and coalescence in the continuum: evolution of states of an infinite
system, Discrete Contin. Dyn. S. 40, 725--752 (2020)

\bibitem{LeJ} Le Jan, Y., Raimond, O.: Flows, coalescernce and
noise. Ann. Probab. 32, 1247--1315 (2004)



\bibitem{Mu}   Murrell, D. J.,  Dieckmann,  U.,  Law,  R.:  On moment closures for population
dynamics in continuous space,   J. Theoret. Biol.  229, 421--432
(2004)


\bibitem{Neuhauser}  Neuhauser, C.: Mathematical challenges in spatial ecology,
Notices of AMS  48 (11), 1304--1314  (2001)


\bibitem{Omel} Omelyan, I., Kozitsky, Yu.: Spatially inhomogeneous
population dynamics: beyond the mean field approximation. J. Phys.
A.: Math. Theor. 52, 305601 (18pp) (2019)
\bibitem{KP} Pilorz, K.: A kinetic equation for repulsive coalescing
random jumps in continuum. Ann. Univ. Mariae Curie-Sk{\l}odowska
Sect. A 70, 47--74 (2016)

\bibitem{Prig} Prigigine, I., Herman, R.: Kinetic Theory of
Vehicular Transport, Elsevier, NY (1971)

\bibitem{Ruelle} Ruelle, D.: Superstable interactions in classical statistical mechanics, Comm. Math. Phys. 18,
127--159  (1970)

\bibitem{Uchiyama} Uchiyama, K.: Derivation of the Boltzmann equation
from particle dynamics, Hiroshima Math. J. 18, 245--297 (1988)



\end{thebibliography}
